\documentclass[10pt]{article}
\usepackage{amsfonts,amsmath,latexsym}
\usepackage{amstext}
\date{}
\newcommand{\R}{\mathbb R} 
 
\newcommand{\N}{\mathbb N} 

\newcommand{\calC}{\mathcal C}
\newcommand{\calD}{\mathcal D}
\newcommand{\calE}{\mathcal E}
\newcommand{\calK}{\mathcal K}
\newcommand{\calI}{\mathcal I}
\newcommand{\calN}{\mathcal N}
\newcommand{\calL}{\mathcal L}
\newcommand{\calM}{\mathcal M}
\newcommand{\calF}{\mathcal F}
\newcommand{\calP}{\mathcal P}
\newcommand{\calO}{\mathcal O}
\newcommand{\calS}{\mathcal S}
\newcommand{\calT}{\mathcal T}
\newcommand{\calU}{\mathcal U}

\newcommand{\calW}{\mathcal W}
\newcommand{\calQ}{\mathcal Q}
\newcommand{\Sig}{\Sigma}
\newcommand{\Sigb}{\overline{\Sigma}}  
\newcommand{\del}{\partial}
\newcommand{\e}{\varepsilon}

\def\ds{\displaystyle}

\newtheorem{theorem}{Theorem} 
\newtheorem{lemma}{Lemma} 
\newtheorem{proposition}{Proposition} 
\newtheorem{corollary}{Corollary} 
\newtheorem{definition}{Definition} 
\newtheorem{remark}{Remark} 

\begin{document} 

\title{The conformal theory of Alexandrov embedded constant mean 
curvature surfaces in $\R^3$}

\author{Rafe Mazzeo\thanks{Supported by the NSF under grant DMS-9971975 and 
at MSRI by NSF grant DMS-9701755} 
\\ Stanford University \and 
Frank Pacard\thanks{Supported at MSRI by the NSF under grant DMS-9701755}
\\ Universit\'e de Paris XII \and  
Daniel Pollack\thanks{Supported by the NSF under grant DMS-9704515
and at MSRI by NSF grant DMS-9701755} 
\\ University of Washington}

\maketitle 

\begin{abstract}
We first prove a general gluing theorem which creates new nondegenerate 
constant mean curvature surfaces by attaching half Delaunay surfaces 
with small necksize to arbitrary points of any nondegenerate CMC 
surface. The proof uses the method of Cauchy data matching from
\cite{MP}, cf.\ also \cite{MPP}. In the second part of this paper, 
we develop the consequences of this result and (at least partially) 
characterize the image of the map which associates to each complete, 
Alexandrov-embedded CMC surface with finite topology its associated 
conformal structure, which is a compact Riemann surface with a finite number 
of punctures. In particular, we show that this `forgetful' map is surjective
when the genus is zero. This proves in particular that the CMC moduli space 
has a complicated topological structure. These latter results are
closely related to recent work of Kusner \cite{Ku}. 
\end{abstract}

\section{Introduction}

In this paper we consider the class of surfaces $\Sigma\subset \R^3$ which
have constant mean curvature (CMC) equal to $1$, are complete and of finite
topology. We shall assume that these surfaces are properly immersed, and in
fact  satisfy the stronger condition that they are {\it Alexandrov embedded}.
This last condition means that 
the immersion $\varphi: \Sigma \hookrightarrow \R^3$ extends to a 
proper immersion $Y \hookrightarrow \R^3$ of a three manifold $Y$
whose boundary is the surface $\Sig$, $\partial Y=\Sig$. 
We remark that, except for the 
sphere, no compact CMC surface, including Wente tori,  ever satisfy this 
condition. The space of all surfaces of this type,  
of genus $g$ with $k$ ends, will be denoted $\calM_{g,k}$. 
To be definite, we do not identify 
elements of this space which differ by Euclidean motions. 

\medskip

The simplest examples of surfaces of this type are the rotationally  invariant 
Delaunay surfaces. Up to Euclidean motions, these are parametrized  by a 
single `necksize' parameter $\tau$, and will be denoted $D_\tau$. We  discuss 
these at greater length later. A remarkable theorem of Meeks \cite{Mee} 
states that each end of an Alexandrov embedded surface is cylindrically 
bounded, and using this Korevaar, Kusner and Solomon \cite{KKS} 
proved that each end is in fact strongly  convergent to some Delaunay 
surface, and thus each end has an associated asymptotic necksize  parameter. 
This strong control on the asymptotic geometry of  these surfaces makes it 
possible to understand the rudimentary structure of these moduli spaces, 
and in 
\cite{KMP} it is proved that  $\calM_{g,k}$ is always a locally real 
analytic space of virtual dimension $3k$. Somewhat remarkably, this 
dimension only depends on the number of ends and not the genus. 
Furthermore, if $\Sig \in \calM_{g,k}$ satisfies a certain analytic 
nondegeneracy condition, which we explain below, then the moduli space 
is a smooth real analytic manifold of dimension $3k$ in a neighbourhood of the 
point $\Sig$. 

\medskip

To describe this nondegeneracy condition, recall that the Jacobi operator
$L_\Sig$ of a CMC surface $\Sigma$ is the linearization of the mean 
curvature operator at $\Sigma$. Surfaces which are $\calC^2$ close to 
$\Sigma$ may be parametrized as normal graphs, i.e.\ as the images
of the map
\[
\Sig \ni x \longmapsto x + w(x)\nu(x)
\]
where $\nu$ is the unit normal to $\Sig$ at $x$ and $w$ is the
(scalar) displacement function. With respect to this representation, $L_\Sigma$ 
takes the simple form $\Delta_\Sigma + |A_\Sigma|^2$. Functions in the 
nullspace of $L_\Sigma$ are called Jacobi fields.  
\begin{definition}
The surface $\Sigma \in {\mathcal M}_{g,k}$ is said to be nondegenerate if 
it has no nontrivial Jacobi fields which decay at all ends of $\Sig$.  
\label{de:nond}
\end{definition}

\medskip

It is natural to investigate 
the structure of the moduli spaces $\calM_{g,k}$ in greater 
detail. Two basic questions are~: first, for which values of $g$ and $k$ is 
$\calM_{g,k}$ nonempty, and second, when does $\calM_{g,k}$ contain a 
nondegenerate element? We now discuss briefly what is known 
about these questions.

\medskip

That $\calM_{g,1}$ is empty for all $g$ is due to Meeks \cite{Mee}.
Also, $\calM_{g,2}$ is empty unless $g=0$, and in this case this 
space contains only 
the Delaunay surfaces. As we indicate  later, Delaunay surfaces are always 
nondegenerate. (As a check on dimensions, note that there is a one 
dimensional family of Delaunay surfaces with fixed  axis which satisfy an 
additional `positioning' normalization determining their translational 
location along this axis. 
The group of Euclidean motions is six dimensional, but the 
rotations about the fixed axis act trivially on these surfaces. 
Thus there is a $6 
= 3 \cdot 2$ dimensional family of  surfaces with $2$ ends, as predicted.)  To 
proceed further it is necessary to use transcendental methods, in particular, 
analytic gluing constructions of these surfaces. The first such surfaces with  
$k\geq 2$ were constructed by Kapouleas \cite{Ka}, but while his method 
gives elements in $\calM_{g,k}$ for infinitely many values of $g$ and $k$,  
it provides little geometric or analytical control on the surfaces themselves, 
and more specifically it seems hard  to determine whether the surfaces he 
constructs satisfy the  nondegeneracy condition.  More recently, the first  
two authors \cite{MP} introduced a new method to handle these gluing 
constructions which involves matching Cauchy data across the gluing 
interfaces. This technique has many advantages over previous  methods~: it 
involves considerably fewer technicalities (and thus  makes it possible to 
approach more complicated geometric problems of this type), it provides very 
good control of the  resulting surfaces, and most importantly, one may often 
prove that the surfaces obtained in this way are nondegenerate. The main 
result  of \cite{MP} is that if $M$ is a complete {\it minimal} surface of 
finite  total curvature in $\R^3$, of genus $g$ with $k$ (asymptotically 
catenoidal, but not planar) ends which is Alexandrov-embedded and  satisfies 
an analogous nondegeneracy condition, then it is possible to obtain a 
nondegenerate CMC surface $\Sig$ of genus $g$ with  $k$ ends by gluing 
half-Delaunay surfaces onto the boundaries of a truncation of $M$. This 
construction yields the existence of nondegenerate elements in 
$\calM_{g,k}$ for infinitely many values of $g$ and $k$.  This result
is a basic ingredient in the nice result of Grosse-Brauckman, Kusner 
and Sullivan \cite{GKS}, who prove that modulo rigid motions, 
$\calM_{0,3}$ is a three-dimensional ball; in particular this space is 
connected. 
(Unfortunately, their method does not prove that every element in this space 
is nondegenerate,  although this is quite likely true.) 
As a further development of this new gluing method, a connected sum 
construction for compact CMC surfaces with boundary
is given in \cite{MPP}. The forthcoming paper \cite{MPPR} extends this 
construction to include connected sums of nondegenerate surfaces in 
$\calM_{g,k}$ and moreover shows that for every $g \geq 0$ and $k \geq 3$,
$\calM_{g,k}$ is nonempty and contains nondegenerate elements. 

\medskip

Beyond these theorems, very little else is known about these 
moduli spaces. It is of interest to determine even the most basic 
properties of  their
topological structure. For example, it is even unknown whether these 
spaces are ever disconnected. It is also unknown whether the 
natural Lagrangian structure \cite{KMP} on these moduli spaces
can be used in any significant way. 

\medskip

The present paper has two main parts which we describe now in turn. 
In the first, which occupies the bulk of the paper, we establish a new 
gluing theorem~: 
\begin{theorem} 
Let $\Sigma \in {\mathcal M}_{g,k}$ be nondegenerate. Then for any point 
$p \in \Sig$ there is a one-parameter family of nondegenerate CMC surfaces 
$\Sigma_\tau(p) \in {\mathcal M}_{g,k+1}$ obtained by gluing a 
half-Delaunay surface $D_\tau$, with $\tau$ sufficiently small, to $\Sigma$ 
at $p$. 
\label{th:3}
\end{theorem}
The geometry of $\Sig$ away from the point $p$ is perturbed very
little in this construction, and in fact as $\tau \to 0$, $\Sig_\tau(p)$
converges, on compact subsets of ${\mathbb R}^3 -\{p\}$,
 to the union of the initial surface $\Sig$ and an infinite
family of mutually tangent spheres of radius $2$ arranged along 
a ray normal to $\Sig$ at $p$.  We actually prove a slightly stronger
theorem~: we show that if $\Sig \in \calM_{g,k}$ is nondegenerate and 
$p \in \Sig$, then there are two distinct one-parameter families of 
nondegenerate CMC surfaces $\Sigma_\tau^\pm(p)$ with $k+1$ ends; these 
are obtained by gluing half of either an embedded Delaunay surface 
(an unduloid) or an immersed Delaunay surface (a nodoid) with very 
small neck. The surfaces $\Sig_\tau^-(p)$ obtained by gluing 
a nodoid on to $\Sig$ are no longer Alexandrov embedded, but we see that
they behave very much like Alexandrov embedded CMC surfaces.
(As an aside, this property of nodoids behaving like Alexandrov embedded 
CMC surfaces is only true when the necksize is small; the
paper \cite{MP2} shows that this fails rather dramatically when
nodoids with large necksizes are considered.)
We remark also that this construction was prefigured and motivated by 
some numerical and computer graphic studies carried out (and brought
to our attention) several years ago by Grosse-Brauckmann.
There are nice illustrations of these surfaces at
http://www.gang.umass.edu/cmc/.

\medskip

An immediate consequence of this gluing theorem is the 
\begin{corollary}
If $\calM_{g,k}$ contains a nondegenerate element, then for 
any $k' > k$, the moduli space $\calM_{g,k'}$ also contains 
a nondegenerate element. 
\end{corollary}

\medskip

The second part of this paper shows how Theorem~\ref{th:3}
can be used to obtain some information about the global
structure of these CMC moduli spaces. To describe these
results, recall  from \cite{KKS} that a complete Alexandrov embedded 
CMC surface of finite topology is conformally equivalent to the 
complement of a finite number of points in a compact Riemann surface.
Thus  any $\Sig \in \calM_{g,k}$ 
is conformally equivalent to a punctured compact
Riemann surface $\Sigb - \{p_1, \ldots, p_k\}$. 
This suggests that it might be useful to consider the
`forgetful map'
\[
\calF_{g,k} = \calF: \calM_{g,k} \longrightarrow \calT_{g,k}.
\]
By definition, $\calT_{g,k}$ is the Teichm\"uller space of conformal   
structures on a surface of genus $g$ with $k$ punctures, and this 
map is defined by sending $\Sig \in \calM_{g,k}$ to the marked 
conformal structure $[\Sig] \in \calT_{g,k}$ determined by its 
(Alexandrov) embedding. In other words, this map forgets the 
CMC embedding of this surface and only retains its conformal structure. 

To begin, we prove the
\begin{theorem} For each $g$ and $k$, the map $\calF_{g,k}$ is real analytic. 
\label{th:ram}
\end{theorem}
This statement requires clarification. As already indicated, $\calM_{g,k}$ 
is a locally real analytic space. What this means is that for any
$\Sig \in \calM_{g,k}$ there is 
an open finite dimensional real analytic manifold $Y$ in the 
space of all Alexandrov-embedded surfaces near to $\Sig$ (in an appropriate 
topology) which contains a neighbourhood $\calU$ of $\Sig$ in 
$\calM_{g,k}$, such that $\calU$ is closed in $Y$ and is given
as the zero set of a real analytic function in $Y$. When
$\Sig$ is nondegenerate, then $Y$ can be taken as the neighbourhood
$\calU$ in $\calM_{g,k}$ itself. To say that $\calF$ is analytic
on $\calM_{g,k}$ means that it has a real analytic extension to
all such real analytic manifolds $Y$. 

\medskip

One consequence of the real analyticity of these spaces and maps
is that $\calM_{g,k}$ is stratified by open real analytic manifolds
$S_j$ such that on each of these strata $\calF$ has constant rank. 

\medskip

We are not claiming that $\calF$ is a proper mapping, and indeed,
it is clear (from any of the gluing constructions) that this is false.  
In another paper in this volume \cite{Ku}, Kusner shows that
there is a suitable modification of the forgetful map which is proper.

\medskip

Recall next that the dimension of $\calM_{g,k}$ around nondegenerate
points is $3k$, while on the other hand, $\calT_{g,k}$ is $6g-6+2k$
dimensional. This suggests that $\calF$ might conceivably be
surjective. This is most likely false, in general, but our next
theorems address this issue~:
\begin{theorem} 
Let $g=0$; then for any $k \geq 3$, the map $\calF_{0,k}$ is surjective.
\label{th:stabgen0}
\end{theorem}

\noindent Kusner has obtained a different proof of this result \cite{Ku} 
using his theorem about the properness of $\calF$.

\begin{theorem}
Fix $g \geq 1$ and suppose that $\calM_{g,k_0}$ contains a 
nondegenerate element for some $k_0 \geq 3$. Then for any $k \geq k_0$,
the image of $\calF_{g,k}$ contains an analytic submanifold of
codimension $d_{g,k}$ which is uniformly bounded as $k \to \infty$.
In other words, the codimension of the image of $\calF_{g,k}$ is
bounded as $k \to \infty$ for each $g$. 
\label{th:stabim}
\end{theorem}
To explain the statement of this last theorem, recall the stratified
structures of these spaces and maps. The image $\calI_{g,k}$ of
$\calF_{g,k}$ is again a stratified real analytic space, and we
say that the maximum dimension of any one of these strata in the
image is the dimension of $\calI_{g,k}$. 

\medskip

Proceeding further, we note that the space $\calT_{g,k}$ has rather 
nontrivial topology. In fact, there is a natural mapping 
\[
F':\calT_{g,k} \longrightarrow \calC_{g,k}
\]
which carries the marked Riemann surface $[(\Sigb;p_1, \ldots, p_k)]$
to the associated element in the configuration space of $k$ distinct
ordered points on a surface of genus $g$. The fundamental group and
cohomology ring of this configuration space have been studied intensively,
see \cite{Bi}, and are known to be rather complicated. It is (barely) 
conceivable that the image of $F'\circ \calF_{g,k} $ does not see 
any of this topology, but this is not the case. 
\begin{theorem}
When $g=0$ and $k \geq 3$, the map 
\[
(F' \circ \calF_{g,k})_*: 
\pi_1(\calM_{0,k}) \longrightarrow \pi_1(\calC(0,k))
\]
is an epimorphism.
If $g \geq 1$ and $\calM_{g,k_0}$ contains a nondegenerate element for 
some $k_0 \geq 3$, then for any $k \geq k_0$, the image of the fundamental 
group $\pi_1(\calM_{g,k})$ under the homomorphism $(F'' \circ \calF_{g,k})_*$ 
contains a finitely generated group with an increasing number of
generators as $k\to \infty$.
\end{theorem}

\medskip

We refer to the final section of this paper for a more detailed statement.

\medskip

The end-to-end gluing construction of Ratzkin \cite{Rat}, 
cf.\ also \cite{MPPR}, also implies the topological nontriviality
of $\calM_{g,k}$, but gives less information than is obtained here.

\medskip

These results together constitute the first and simplest steps of a 
more detailed investigation of the topology of the moduli spaces  
$\calM_{g,k}$.

\medskip

The outline of the rest of the paper is as follows. Sections 2 through 4
contain the proof of the gluing result, Theorem~\ref{th:3}. More specifically,
\S 2 contains the analysis of CMC deformations of half-Delaunay
surfaces and \S 3 contains the analysis of CMC deformations of
$\Sig - D$, where $\Sig$ is any element in $\calM_{g,k}$
and $D$ is a small geometric disk in $\Sig$. These results are then
combined in \S 4, where it is shown how to perform the Cauchy
data matching across the gluing interface, and hence how to produce
a new CMC surface with $k+1$ ends. Finally, in \S 5 we undertake the
analysis of the forgetful map and develop its properties and
derive the theorems stated above concerning the image of
$\calM_{g,k}$. 

\medskip

The authors wish to thank Matthias Weber and Rob Kusner for 
a number of useful conversations. The first author also wishes to
thank Ralph Cohen and Gunnar Carlsson for setting him straight
on the topology of configuration spaces.

\section{The geometry and analysis of Delaunay surfaces}

In this section we recall the family of Delaunay surfaces $D_\tau$,  and 
review some of their basic geometric and analytic properties. We focus 
particularly on their behavior in the singular limit as $\tau \to 0$.  

\subsection{Definition and first properties}

To find all CMC surfaces in $\R^3$ which are rotationally invariant  about an 
axis, one is quickly led to the ODE which the generating curve  for such 
surfaces must satisfy. This ODE has a first integral, and from this one 
can see 
that all of its solutions are periodic. Thus  we obtain a family of periodic, 
rotationally invariant surfaces  $D_\tau$. With a particular normalization 
described in the next paragraph, the parameter $\tau$ lies in $(-
\infty,0)\cup(0,1]$;  when $\tau > 0$,  the surface $D_\tau$ is embedded, 
while when $\tau < 0$, it is  only immersed. There is a geometric limit, 
as $\tau \to 0^{\pm}$, consisting of an infinite arrangement of spheres, 
each tangent to 
the next, arranged along an axis.  We shall not describe this 
more carefully,  but instead refer to \cite{MP} for details of 
the preceding  discussion as well as the material in the remainder of this 
section. 

\medskip

The parametrization given by this cylindrical coordinate description  
is not very tractable analytically, but fortunately it turns out 
that there is a parametrization which is much easier to use. 
This is the isothermal parametrization given by
\begin{equation}
X_\tau: \R \times S^1 \ni (s,\theta) \longmapsto  
\frac12\,\left(\tau\,e^{\sigma(s)}
\, \cos\theta,\tau\,e^{\sigma(s)}\,  \sin\theta,\kappa (s)\right),
\label{eq:2.1}
\end{equation}
where the functions $\sigma $ and $\kappa $ are described as follows
(cf. \S 3 of \cite{MP}). For any $\tau \in (0,1]$, the function $\sigma $ 
is defined to be the unique smooth nonconstant solution of the ODE
\[
(\del_s\sigma)^2 + \tau^2 \cosh^2 \sigma =1, \qquad \del_s \sigma(0) = 0,  
\qquad \sigma(0) = -\operatorname{arccosh} 1/\tau ,
\]
while, for any $\tau \in (-\infty,0)$, 
the function $\sigma $ is defined to be the 
unique smooth nonconstant solution of the ODE
\[
(\del_s\sigma)^2 + \tau^2 \sinh^2 \sigma =1, \qquad \del_s \sigma(0) = 0,  
\qquad \sigma(0) = \operatorname{arcsinh} 1/\tau .
\]
Again, the definition of $\kappa $ differs according to whether $\tau$  is 
positive or negative. If $\tau \in (0,1]$, then we define the 
function $\kappa$ by 
\[
\del_s\kappa =\tau^2\, e^{\sigma}\, \cosh\sigma, \qquad \kappa(0) = 0,
\]
while if $\tau < 0$,  we define the function $\kappa$ by 
\[
\del_s \kappa =  \tau^2 \, e^{\sigma}\,\sinh \sigma, \qquad \kappa(0) = 0.
\]

Observe that when $\tau>0$, $\kappa$ is monotone increasing,  and hence 
$X_\tau$ is an embedding, whereas when $\tau < 0$, this is no longer true and 
the surfaces are only immersed. The Delaunay  surfaces $D_\tau$ are known 
as unduloids and nodoids in these two  cases, respectively.  The extreme 
element in the family of unduloids is $D_1$, the cylinder  of radius $1/2$. 
The 
limit of $D_\tau$, either as $\tau \nearrow 0$ or as $\tau \searrow 0$, is an 
infinite union of tangent spheres of  radius $1$, arrayed along a common axis. 

\medskip

As noted above, these surfaces are all periodic. When $\tau > 0$,  there is a 
unique $t_\tau > 0$ such that $\tau \, \cosh t_\tau = 1$,  and if we define
\begin{equation}
s_\tau := \frac12 \,\, \int_0^{t_\tau}\frac{dt}{\sqrt{1-\tau^2\, \cosh^2 t}},
\label{eq:2.2}
\end{equation}
then it is not hard to see that $\sigma $ has period $8\, s_\tau$.   
On the other 
hand, when $\tau < 0$, we define $t_\tau>0$ by the equation $\tau\, \sinh 
t_\tau = -1$ and define
\begin{equation}
s_\tau := \frac12 \,\,\int_0^{t_\tau}\frac{dt}{\sqrt{1-\tau^2\, \sinh^2 t}} ,
\label{eq:2.3}
\end{equation}
once again $\sigma $ is periodic of period $8\, s_\tau$.

\subsection{The singular limit of $D_\tau$ as $\tau \to 0$}
\label{singlim}

We now describe some aspects of the surfaces $D_\tau$ as $\tau \to 0$. In 
this limit, $D_\tau$ converges to a singular `noded' surface,  which is the 
infinite union of spheres of radius $1$ centered at the points  
$(0,0,2 \, k +1)$, 
$k \in {\mathbb Z}$. From (\ref{eq:2.2}) or (\ref{eq:2.3}) we see that 
\[
s_\tau  = - \frac{1}{4} \log \tau^2  + {\mathcal O}(1)  \qquad 
\mbox{as}\quad \tau \to 0.
\]

To study this limit more closely, notice that the family of rescaled  surfaces 
$\tau^{-2}\, D_\tau$ converges to a catenoid of revolution  around the 
$z$-axis. Thus when $\tau$ is small, $D_\tau$ is well approximated by a 
sequence of spheres along the $z$-axis connected by  small catenoidal necks. 
This is central in much of what follows, and so in the remainder of this 
subsection we make this quantitative by recalling the behavior of the 
functions $\sigma$ and $\kappa$ as  $\tau \to 0$.  
The expansions below are not hard to derive, but we  refer to 
\cite{MP} for detailed proofs.

\begin{definition} 
The notation $g = {\calO}_{\calC^\infty}(f)$  means that for any $k \geq 0$ 
there exists $c_k >0$ such that 
\[
|\del^k_s g |\leq c_k \, |f|
\]
on the domains of definition of these functions. 
\end{definition}

To obtain the asymptotics of $\sigma $ as $\tau \to 0$, define 
$r (s)  := \tau \, e^{\sigma(s)}$. 
This function satisfies the ODE 
\[
(\del_s r)^2  = r^2 - \left( \frac{\tau^4}{16} \pm 
\frac{\tau^2}{2}\, r^2 +  r^4 \right),
\]
where the plus or minus sign is chosen according to the sign of $\tau$ ($+$ 
when $\tau >0$ and $-$ when $\tau <0$). Because $\sigma$ attains its  
minimum when $s=0$, $ r (s)\geq  r (0)$  for all $s$. Hence we can write 
$r (s) := r (0) \, \cosh w(s)$ for  some function $w$. This new function 
satisfies 
\[
(\del_s w)^2=1\mp\frac{\tau^2}{2}-  r(0)^2\,(\cosh^2 w+1).
\]
Now simply from the definition of $\sigma(0)$ we have 
\[
r (0) = \frac{\tau^2}{4} + {\mathcal O} (\tau^4) \qquad  \mbox{as}\quad 
\tau \to 0,
\]
and therefore
\[
w(s)=\left(1\mp\frac{\tau^2}{2}\right)^{1/2}\, s+{\calO}_{\calC^\infty}
(\tau^4\,\cosh^2 s).
\]
Since $s_\tau \sim -\frac14 \log \tau^2$, $\tau^4 \cosh^2 s \leq c$  when 
$|s| \leq 4s_\tau$, and so this estimate is only interesting  when 
$s \in [- 4\, s_\tau,4\, s_\tau]$. 
However, by periodicity  we obtain an estimate for $w$ for all $s \in 
\R$. From here it  follows that
\[
\kappa (s)=  \pm \frac{\tau^2}{2}\, s+\calO_{\calC^\infty}(\tau^4 \,\cosh^2 s),
\]
according to the sign of $\tau$ (with $+$ when $\tau >0$ and $-$ with $\tau 
<0$) but this is only of interest when $s \in [-2\,s_\tau,2\, s_\tau]$.  
It is possible to refine this argument to obtain a more precise expansion for 
$\kappa$ in the larger interval $[-4\, s_\tau,4\, s_\tau]$,  
but we omit this since it will not be needed. 

\medskip

The image of $X_\tau$ restricted to $(0,2\, s_\tau]\times S^1$   may also be 
written as the graph of functions $U_\tau$, which is defined over an annulus 
in the $x \, y$-plane.  (The image of $X_\tau$  restricted to $[-2s_\tau,0)\times 
S^1$ is the graph of $-U_\tau$.)  To analyze $U_\tau$ as $\tau \to 0$, we
use the function $r =\frac{\tau}{2} \, e^{\sigma(s)}$, as above and set 
\[
r_\tau  :=\frac{\tau}{2} \, e^{\sigma (-s_\tau)}.
\]
For the moment, observe that $r_\tau \sim \tau^{3/2}$ as $\tau$ tends to 
$0$. Expanding  $\kappa$ in terms of $r$ near $r=r_\tau$, we find that when 
$r_\tau/2 \leq r \leq 2r_\tau$,  
\begin{equation}
\left| (r\,\nabla)^k\,\left(U_\tau(r,\theta) \mp \frac{\tau^2}{4}\,\log\left( 
\frac{8r}{\tau^2}\right)\right) \right| \leq c_k \, \tau^3 ,  \qquad  
\forall k \geq 0,
\label{eq:2.4}
\end{equation}
according to the sign of $\tau$ (with $-$ when $\tau >0$ and $+$ with $\tau 
<0$).  Similar estimates are valid on a larger set, but this behavior near $r 
=r_\tau$ will be the most crucial later.

\subsection{The Jacobi operator on a Delaunay surface}

If $\Sigma$ is a CMC surface, then any surface which is $\calC^2$ close  to it 
may be represented as a normal graph 
\[
\Sigma_w = \{x + w(x)\nu(x)\quad | \quad  x \in \Sigma\},
\]
where $\nu$ is the unit normal vector field and $w$ is a (small) scalar 
function. $\Sigma_w$ is itself CMC provided $w$  satisfies a nonlinear 
second order elliptic equation. This equation will be discussed in more detail 
later, but it is well-known that its linearization, $\calL$, 
which is usually called the Jacobi operator, 
is the sum of the Laplace-Beltrami operator on $\Sigma$ 
and the norm squared of the second fundamental form, 
\[
\calL=\Delta_\Sigma + |A_\Sigma|^2.
\]

\medskip

A solution $w$ of the equation $\calL w = 0$ is called a Jacobi field on 
$\Sigma$. Ideally, such a function is the tangent vector of a  one-parameter 
family of CMC deformations of $\Sigma$. This may not be  the case, 
however, and in general, if $w$ is a Jacobi field, then the elements in
any one-parameter family of surfaces $\{\Sigma_{w(\e)}\quad | \quad |\e| < \e_0\}$,
where $w'(0) = w$, only satisfy the constant mean curvature equation
to second order at $\e = 0$. 
Nonetheless, an understanding of the Jacobi fields  and mapping properties of 
the Jacobi operator is fundamental to any account of the deformation theory of 
$\Sigma$. 

\medskip

We denote by $\calL_\tau$ the Jacobi operator associated to the Delaunay  
surface $D_\tau$. In terms of the isothermal parametrization from \S 2.1,
\[
\calL_\tau = \frac{4}{\tau^2e^{2\sigma}}\left(\del_{s}^2 +  \del_{\theta}^2 + 
\tau^2 \, \cosh (2\sigma) \right). 
\]
We analyze instead the simpler operator
\begin{equation}
L_\tau := \del^2_{s} + \del^2_{\theta} + \tau^2 \, \cosh (2\sigma), 
\label{eq:3.1}
\end{equation}
with the factor $4/(\tau^2 e^{2\sigma})$ removed; it is clear that the 
mapping properties of one of these operators implies the corresponding 
properties of the other, and also, their nullspaces are the same.  
Observe that, since $\sigma$ is $8 \, s_\tau$ periodic and even, 
the potential in $L_\tau$ is $4 \, s_\tau$-periodic.

\medskip

We now undertake a thorough analysis of the operators  $\calL_\tau$ and 
$L_\tau$. 

\subsection{Jacobi fields on $D_\tau$}

The uniqueness theorem for 2-ended Alexandrov embedded CMC surfaces \cite{KKS}
implies that the only CMC deformations of an entire Delaunay surface are 
the obvious ones, namely those arising from rigid motions of 
$\R^3$ and changes in the Delaunay 
parameter. The corresponding space of Jacobi fields are either bounded or 
linearly growing along the Delaunay axis. All other Jacobi fields grow 
exponentially in one direction or the other along this axis, and hence do 
not correspond to actual CMC deformations. In the next subsection we 
describe this former class of `geometric'  Jacobi fields, while in 
the one following that we discuss some features of the latter class. While 
these do not lead to global CMC deformations, they do correspond to CMC 
deformations over half of $D_\tau$, and in any case, it is necessary to 
understand them in order to describe the mapping properties of $\calL_\tau$. 

\subsubsection{Geometric Jacobi fields} 

As indicated above, there is a special collection of Jacobi fields on  
$D_\tau$ 
which correspond to explicit geometric deformations of this surface.  This 
family is six-dimensional; there is a three-dimensional space  associated to 
translations in $\R^3$, a two-dimensional space associated to the rotations of 
the Delaunay axis, and a one-dimensional space associated to changing the 
Delaunay parameter. We now describe  six Jacobi fields, which we denote by 
$\Phi_\tau^{j,\pm}$, $j = 0, \pm 1$,  which form a basis for this space.  
The motivation for this notation will be made apparent in subsection
\ref{expjacfields}.

\medskip

We first treat the case where $\tau \neq 1$. 

\medskip

Let us start with $\Phi_\tau^{0,+}$. This Jacobi field corresponds to an 
infinitesimal translation of $D_\tau$ along its axis, and so  is obtained by 
projecting the constant vector field  $(0,0,1)$ (which is the Killing field 
associated to this family of  translations) along the normal vector field 
$N_\tau$ on $D_\tau$. It is  geometrically obvious that $\Phi_\tau^{0,+}$ 
depends only on $s$, and  is periodic in $s$, hence is bounded as $s \to \pm 
\infty$. 

\medskip

Next, the two Jacobi fields $\Phi_\tau^{\pm 1, +}$ correspond to 
translations of $D_\tau$  in the 
two directions orthogonal to its axis. These are calculated by projecting the 
constant vector fields $(1,0,0)$ and $(0,1,0)$ along $N_\tau$, and hence are 
once again bounded in $s$. Moreover,  
\[
\Phi_\tau^{1,+} (s, \theta) = \phi_\tau^{1,+}(s)\cos \theta, \qquad 
\mbox{and}\qquad \Phi_\tau^{-1,+} (s, \theta) = \phi_\tau^{-1,+}(s)\sin 
\theta. 
\]

Continuing, the two Jacobi fields 
\[
\Phi_\tau^{1,-} (s, \theta) = \phi_\tau^{1,-}(s)\cos \theta, \qquad  
\mbox{and}\qquad \Phi_\tau^{-1,-} (s, \theta) = \phi_\tau^{-1,-}(s)\sin \theta 
\]
correspond to infinitesimal rotations of $D_\tau$ about its axis. These are 
given by projecting the Killing fields $(z,0,-x_1)$ and $(0,z, -x_2)$ 
(corresponding to rotations in the space $\R^3$ with coordinates 
$(x_1,x_2,z)$) along $N_\tau$, and hence grow linearly in $s$. 

\medskip

Finally, $\Phi_\tau^{0,-}$ corresponds to the derivative of the one parameter
family $D_\tau$  obtained by varying the Delaunay parameter $\tau$.
Since $D_\tau$ are surfaces of revolution, this Jacobi field depends only 
on $s$. It is linearly growing when $\tau \neq 1$. 

\medskip

The definitions for $\Phi_{\tau}^{0,-}$,  $\Phi_{\tau}^{\pm 1,-}$
and $\Phi_{\tau}^{\pm 1,+}$ when 
$\tau=1$ are exactly the same, but since $D_1$ is a cylinder, the definition 
for $\Phi_1^{0,+}$ above yields $0$.  There is an intrinsic way to 
define the Jacobi fields $\Phi^{0, \pm}_\tau$  for all $\tau \in (0,1]$ by 
regarding the change of Delaunay parameter and translation along the axis as 
playing the role of polar coordinates  in the space of parameters. More 
precisely, let $(\rho,\alpha)$ be the ordinary polar coordinates 
(corresponding to Cartesian coordinates  
$(y_1,y_2) = (\rho\cos\alpha,\rho\sin\alpha)$). Then 
there is a smooth map  from a small ball in $\R^2$ into the two-dimensional 
space of Delaunay  surfaces sharing the same axis, given by 
\[
[0,1)\times [0,2\,\pi) \ni (\rho,\alpha) \longrightarrow D_{1-\rho}+ 
\left(0,0,\frac{\kappa_{1-\rho}(8 \, s_{1-\rho}) }{4 \, \pi}\, \alpha\right) 
\]
Then for any value of $\tau$, the Jacobi fields $\Phi^{0,\pm}_\tau$ have  the 
same span as the variations of this family with respect to the Cartesian 
coordinates. 

\medskip

Further details of these calculations, as well as explicit expressions  
of these geometric Jacobi fields in terms of $\sigma $ and $\kappa $  
can be found in \cite{MP}. 

\subsubsection{Jacobi fields of exponential type and indicial roots} 
\label{expjacfields}

We now fit the special family of Jacobi fields discussed above into  
a broader context.

\medskip

The operator $L_\tau$ is invariant with respect to rotations about  the 
Delaunay axis, and is reduced by the eigendecomposition for the 
cross-sectional operator $\del_\theta^2$. This reduces the  analysis of 
$L_\tau$ to that of the countable family of  operators  
\[
L_{\tau,j}:= \del_{s}^2+(\tau^2\,\cosh(2\sigma)-j^2), \qquad j =0, 1, 2, 
\ldots.
\]
When $j \neq 0$, the operator $L_{\tau,j}$ occurs  
with multiplicity 
two, corresponding to the two eigenfunctions $e^{\pm ij\theta}$.

\medskip

Except in certain exceptional cases described below, for each $j$  
there exists a complex number $\zeta_{\tau,j}$, with $\Re  \zeta_{\tau,j}\geq 0$,
and a basis of solutions $\Psi_\tau^{j,\pm}$ of each of these ordinary differential  
operators such that 
\[
\Psi_\tau^{j,\pm}(s+4s_\tau) = e^{\pm 4 \, \zeta_{\tau,j} \, s_\tau} \, 
\Psi_\tau^{j,\pm}(s).
\]
Thus when $\zeta_{\tau,j}$ is real, $\Psi_\tau^{j,\pm}$  grows exponentially 
as $s \to \pm \infty$ and decays exponentially  as $s \to \mp \infty$. 
It happens that either $\zeta_{\tau,j}$ is real or it is pure imaginary. 
In this case, the corresponding  solutions are oscillatory, hence bounded. 

\medskip

The exceptional cases noted above occur when there is one bounded 
(oscillatory) solution and one solution which grows  linearly. For example, 
we have already seen that this is the case when  $j=0,\pm 1$.

\medskip

One way to prove these facts is to recall that since the potential in 
$L_{\tau,j}$  has period $4s_\tau$, there exists a $2\times 2$ matrix 
$T_{\tau,j}$ such that for any solution $w$ of $L_{\tau,j}w = 0$ on $\R$, 
\[
\left(
\begin{array}{rllll}
w(s + 4 \, s_\tau)\\[3mm]
\del_s w(s+4 \, s_\tau)
\end{array}
\right) = T_{\tau,j} \, \left(
\begin{array}{rllll}
w(s)\\[3mm]
\del_s w(s)
\end{array}
\right).
\]
Using the Wronskian related to the operator $L_{\tau,j}$, it is easy 
to see that the determinant of this matrix is equal to $1$.  Since the matrix 
$T_{\tau,j}$ has real entries, we see that the roots of its characteristic
 polynomial  $\lambda_{\tau,j}^\pm$  satisfy  $\lambda_{\tau,j}^+ \, 
\lambda_{\tau,j}^-  =1$, and  $\lambda_{\tau,j}^+ + \lambda_{\tau ,j}^- 
\in{\mathbb R}$. So we can write them as 
\[
\lambda_{\tau,j}^{\pm} = e^{\pm 4 \, \zeta_{\tau,j}\, s_\tau},
\]
where $ \zeta_{\tau,j}$ is either real or purely immaginary. In the case 
where $\lambda_{\tau, j}^\pm =  1$ or $-1$ it may happen that the matrix 
$T_{\tau,j}$ cannot be diagonalized and this cooresponds to the 
exceptional cases we were mentionning above.  We leave the details of 
checking the statements in the last paragraph 

\begin{definition}
For each $j$ we define the indicial roots of the operator $L_{\tau, j}$  to be 
the pair of numbers $\pm \gamma_{\tau,j}$ where 
\[
\gamma_{\tau,j} : = \mbox{Re\,} \zeta_{\tau,j} \geq 0.
\]
Thus 
\[
\Gamma_{\tau} := \left\{\pm \gamma_{\tau,j} \quad | \quad  j \in {\mathbb N} \right\}
\]
is the set of all indical roots of the operator $L_{\tau}$. 
\end{definition}

It is not necessary to express these indicial roots exactly 
(and indeed,  it is probably impossible to do so), but the 
following estimates will suffice for our nefarious purposes.
\begin{proposition}
The indicial roots of $L_\tau$ satisfy the following properties~:
\begin{itemize}
\item[(i)] For any $\tau \in (-\infty,0)\cup (0,1]$,  $\gamma_{\tau,0} 
=\gamma_{\tau,1} =0$. 
\item[(ii)] When $j \geq 2$ and $\tau \in (-\sqrt{j^2-2},0) \cup (0,1]$, 
$\gamma_{\tau,j} > 0$ and $L_{\tau,j}$ satisfies the maximum principle.
\item[(iii)] For any $j \geq 2$, $\lim_{\tau\rightarrow 0} 
\gamma_{\tau,j}  = j$.
\end{itemize}
\label{pr:3.1}
\end{proposition}
{\bf Proof:} Property (i) follows from the remarks above. Notice that because 
of the existence of linearly growing solutions,  $\zeta_{\tau,j}$, $j=0,1$, is 
identically zero for all $\tau$, instead  of just having real part zero. 

\medskip

To prove property (ii), it suffices to show that when $\tau$ is in the  
stated range, the potential in $L_{\tau,j}$ is strictly negative  and so this 
operator satisfies the maximum principle. Hence it  cannot have bounded 
solutions and therefore its indicial roots must  be real and nonzero. 
To show this, first assume that $\tau > 0$. Then
\[
\tau^2\, \cosh (2\,\sigma)-j^2 = 2 \,\tau^2 \,\cosh^2\sigma -  \tau^2 -j^2 
= 2 - \tau^2 - j^2 - (\del_s \sigma)^2,
\]
which is strictly negative when $j \geq 2$. On the other hand, 
when $\tau < 0$, 
\[
\tau^2\, \cosh (2\,\sigma)-j^2 = 2 \, \tau^2 \,\sinh^2 \sigma +  
\tau^2 -j^2  =  2 + \tau^2 - j^2 - (\del_s \sigma)^2,
\]
which is strictly negative when $- \sqrt{j^2-2} < \tau < 0$. 

\medskip

Finally, Property (iii) reflects the fact that as $\tau$ tends to $0$,  
the potential in $L_\tau$ is arbitrarily close to $0$ on sets which  
are arbitrary large. To be 
more precise, according to Proposition~13 of  \cite{MP}, for any $\eta 
>0$ there exist numbers $s_\eta >0$ and  $\tau_0 >0$ such that when $| \tau  
| \in (0,  \tau_0)$, 
\[
\tau^2 \, \cosh (2\sigma) \leq \eta ,
\]
when $s \in [s_\eta, 4 \, s_\tau - s_\eta]$. 
It follows from this that the indicial 
roots of $L_{\tau,j}$ must converge to the indicial roots  of the operator 
$\del_s^2 -j^2$ as $\tau \to 0$, see  Proposition~20 in 
\cite{MP} for details.\hfill $\Box$

\medskip

This Proposition proves that there exists a number $\tau_*\leq-\sqrt{2}$ 
such that 
\begin{equation}
j \geq 2 \ \mbox{and}\ \tau \in (\tau_*, 0)\cup (0, 1] \Longrightarrow
\gamma_{\tau,j} >0, \ \mbox{and $L_{\tau,j}$ satisfies the maximum 
principle.}
\label{eq:deftst}
\end{equation}

\subsection{Mapping properties of $L_\tau$}

We shall call the image by $X_\tau$ of $[s_0,\infty) \times S^1$, for any $s_0 
\in \R$, a half-Delaunay surface, and sometimes denote it  by 
$D_\tau^+(s_0)$. Our goal in the remainder of this section
is to study CMC perturbations  of half-Delaunay surfaces with prescribed
boundary values.
To do this we require rather precise knowledge of the mapping 
properties of the Jacobi operators, and their inverses, on these surfaces. 
In fact, we need to 
understand these mapping properties uniformly as $\tau \to 0$. In this 
subsection we give careful statements of  the results we need. Results of this 
type (for fixed $\tau$) were  originally proved in \cite{MPU} for 
$L_\tau$ acting on weighted Sobolev spaces; these were reformulated for the 
more convenient  family of weighted H\"older spaces in \cite{MP}, and 
the issue of uniformity in the singular limit was also addressed there. 

\medskip

Let us begin by defining these weighted H\"older spaces on $D_\tau^+ (s_0)$.
\begin{definition} 
Let $r \in {\mathbb N}$, $\alpha \in (0,1)$, and $\mu \in \R$. Then the 
function space $\calE^{r,\alpha}_\mu([s_0,\infty) \times S^1)$ consists of  
those functions $w \in \calC^{r,\alpha}_{loc}([s_0,\infty) \times S^1)$  such 
that
\[
\|w\|_{\calE^{r,\alpha}_\mu}:= \sup_{s \geq s_0} e^{-\mu s}\, \| 
w\|_{\calC^{r,\alpha}([s,s + 1]\times S^1)} < \infty. 
\] 
Here  $\| \cdot \|_{\calC^{r,\alpha}([s,s+1]\times S^1)}$ is the  
usual H\"older norm on $[s,s +1] \times S^1$. 
Also, $[\calE^{r,\alpha}_\mu ([s_0,\infty)\times 
S^1)]_0$ is the subspace of functions vanishing at $s=s_0$. 
\label{de:4.1}
\end{definition} 
Observe that the function $s \mapsto e^{\mu s}$ is in $\calE^{r,\alpha}_\mu  
(\R^+ \times S^1)$.

\medskip

It is clear that 
\begin{equation}
L_\tau :[{\mathcal E}^{2,\alpha}_\mu ([s_0, +\infty) \times S^1)]_0  
\longrightarrow {\mathcal E}^{0,\alpha}_\mu ([s_0, +\infty) \times S^1)
\label{eq:blt}
\end{equation}
is bounded for all $\mu \in \R$. However, it is not Fredholm for every weight. 
In fact, the existence of a solution of $L_\tau w = 0$  which grows or decays 
exactly like $e^{\pm \gamma_{\tau,j} s}$ can be used to show that this 
mapping does not have closed range when $\mu = \pm \gamma_{\tau,j}$ for 
any $j$. However, it is proved in \cite{MP}, cf. also \cite{MPU}, 
that (\ref{eq:blt}) is Fredholm when $\mu \notin \Gamma_{\tau}$. 

\medskip

We first consider what happens for $\tau$ fixed. For simplicity, denote by 
$L_\tau(\mu)$ the operator in (\ref{eq:blt}); hence this notation  
indicates the 
weighted space on which we are letting $L_\tau$ act. A basic observation is 
that $L_\tau(\mu)$ and $L_\tau(-\mu)$ are essentially dual to one another. 
(Of course, since we are working on H\"older spaces, this is not quite
true. The analogous statement for weighted $L^2$ spaces is true, however,
and this may be used to deduce the following statements.)
An important consequence of this is that when $\mu \notin 
\Gamma_\tau$, then $L_\tau(\mu)$  is surjective if and only if 
$L_\tau(-\mu)$ is 
injective; furthermore, if this is the case, then the dimension of the 
kernel of $L_\tau(\mu)$ is equal to the dimension of the cokernel of 
$L_\tau(-\mu)$.   Summarizing all of this, the precise result is
\begin{proposition} 
Assume that $\mu\in(\gamma_{\tau,j},\gamma_{\tau,j+1})$  for some $j \in 
{\mathbb N}$, then $L_{\tau}(\mu)$ is surjective and has a kernel of 
dimension $2 \, j+1$. Thus when $\mu \in (- \gamma_{\tau, j+1}, - 
\gamma_{\tau, j}) $  for some $j \in {\mathbb N}$, then $L_{\tau}(\mu)$ is 
injective, and has a cokernel of dimension $2 \, j+1$.   
\label{pr:4.1}
\end{proposition} 

Notice that the first of these assertions follows from the discussion  
of Jacobi 
fields in the last subsection. In fact, we need all 
$\Psi_\tau^{+,j} \, e^{\pm ij 
\theta}$ and $\Psi_\tau^{-,j} \, e^{\pm ij \theta}$ to lie in 
$\calE_\mu^{r,\alpha} ([s_0,\infty))$ in order to take a linear combination of 
these which vanishes at $s_0$. 

\medskip

We shall henceforth assume that $\tau \in (\tau_*,0)\cup (0,1]$ and $\mu \in 
(-\gamma_{\tau,2},-\gamma_{\tau,1}) =  (-\gamma_{\tau,2} ,0)$. Although 
the  previous Proposition states that $L_\tau(\mu)$ is not surjective in this 
case, there is a way to modify this mapping so as to be surjective by
augmenting the domain. In fact, from the comments in the  last paragraph, 
it is clear that if 
we augment the domain  of $L_\tau(\mu)$, for $\mu$ in this (negative) range, 
with a `deficiency space' 
\[ 
\calW_\tau := \mbox{Span}\, \{ \Psi^{j,\pm}_\tau(s) e^{ij\theta} \quad | \quad
j=0, \pm 1\},
\] 
which is simply the span of the geometric Jacobi fields, then  
the nullspace of 
$L_\tau(-\mu)$ is contained in this extended domain, i.e.  
\[ 
\mbox{Ker}\, L_{\tau}(-\mu) \subset [{\mathcal E}^{2,\alpha}_{\mu}([s_0,  
\infty)\times S^1) \oplus \calW_\tau]_0.
\] 
The deficiency subspace $\calW_\tau$ is $6$-dimensional, while from  
Proposition~\ref{pr:4.1} again, $\dim\mbox{Ker}\,L_\tau(-\mu) = 3$.  Thus 
projecting this nullspace onto $\calW_\tau$ defines a $3$-dimensional  
subspace ${\cal N}_\tau \subset \calW_\tau$. Choosing any complement 
${\cal K}_\tau$, so that 
\[ 
\calW_\tau = \calN_\tau \oplus \calK_\tau, 
\] 
then 
\begin{equation}
L_\tau: [\calE^{2,\alpha}_{\mu}([s_0,\infty)\times S^1) \oplus \calK_\tau]_0 
\longrightarrow \calE^{0,\alpha}_{\mu}([s_0,\infty)\times S^1)
\label{eq:auglm}
\end{equation}
is injective. With slightly more work, as explained in \cite{MP},  
\cite{MPU}, it is possible to show that  it is also surjective. Altogether, 
we have 
\begin{proposition}  
For $\tau \in (\tau_*, 0) \cup (0,1]$, $\mu \in (-\gamma_{\tau,2}, -
\gamma_{\tau,1})$ (which is the same as $(-\gamma_{\tau,2},0)$ since
$\gamma_{\tau,1} = 0$), and any choice of complement 
$\calK_\tau$ to  $\calN_\tau$ in $\calW_\tau$, then the mapping 
(\ref{eq:auglm}) is an isomorphism. 
\label{pr:4.2} 
\end{proposition} 

It is more convenient for us to rephrase this result using the fact  
that the space of boundary traces $\{w(s_0,\cdot) \quad | \quad w \in \calK_\tau\}$  
lies in the span of $e^{ij\theta}$, $j=0, \pm 1$~:
\begin{proposition} 
Fix $\tau \in (\tau_*, 0) \cup (0,1]$ and $\mu\in(-\gamma_{\tau,2},  -
\gamma_{\tau,1})= (-\gamma_{\tau,2},0)$. Then for any $s_0 \in \R$, there 
exists a bounded mapping 
\[
G_{\tau,s_0}:\calE^{0,\alpha}_{\mu}([s_0,\infty)\times S^1)  \longrightarrow 
\calE^{2,\alpha}_{\mu}([s_0 ,\infty)\times S^1)
\]
such that for any $f \in \calE^{0,\alpha}_{\mu}([s_0,\infty)\times S^1)$,  
the function $w = G_{\tau,s_0}(f)$ satisfies
\[
\left\{ 
\begin{array}{rlll} 
L_\tau  w & = & f \quad  & \mbox{in}\quad (s_0,\infty) \times S^1\\[2mm] 
w & \in   & \mbox{\rm Span}\,\{e^{-i\theta},1, e^{i\theta}\} \quad & 
\mbox{on}\quad \{s_0\} \times S^1.
\end{array} 
\right. 
\]
Moreover the norm of $G_{\tau, s_0}$ is bounded independently of $s_0 \in 
{\R}$.
\label{pr:4.3}
\end{proposition} 

As explained earlier, we must also understand the behaviour of this inverse 
as $\tau \to 0$. To explain the statement of this next result, recall  that, 
for all $j \geq 2$, $\gamma_{\tau,j} \to j$ as $\tau \to 0$, so that any fixed 
$\mu \in (-2,-1)$ is eventually in $(-\gamma_{\tau,2}, -\gamma_{\tau,1}) 
= (-\gamma_{\tau,2}, 0)$ provided $\tau$ is chosen small enough. 
\begin{proposition} 
Fix any $\mu \in (-2,-1)$. Then, there exists a number $\tau_0 >0$ such that 
the norm of $G_{\tau, s_0}$ is uniformly bounded, independently of $s_0 \in 
\R$ and $| \tau  | \in (0, \tau_0)$. 
\label{pr:4.4} 
\end{proposition}
The precise range of $\mu$ between $-2$ and $-1$  is very important here; 
this result fails when, for example,  $\mu \in (-1,0)$. This is  
Proposition~21 in \cite{MP}, and the proof can be found there.

\medskip

To conclude, we need some facts about the Poisson operator. However,  
we only require these in the limit as $\tau \to 0$, and so it suffices 
to study the Poisson operator for the very simple operator 
\[
\Delta_0 : = \del_s^2 + \del^2_\theta
\] 
on $[0, \infty)\times S^{1}$. 
\begin{lemma}
For any $h \in  \calC^{2,\alpha}(S^{1})$ such that 
\[
\int_{S^1}h(\theta)\, d\theta = \int_{S^1} h(\theta)e^{\pm i\theta}\,d\theta 
= 0,
\]
(i.e. the function $h$ is orthogonal to $1$ and $e^{\pm i\theta}$ in the $L^2$ 
sense on $S^1$) there exists a unique function $w = \calP(h) \in  
\calE^{2,\alpha}_{-2}([0, +\infty)\times S^{1})$ such that $\Delta_0 w = 0$ 
on $(0,+\infty)\times S^{1}$,  $w(0,\theta) = h(\theta)$, and moreover, $\| w 
\|_{\calE^{2,\alpha}_{-2}} \leq c \, \|h\|_{\calC^{2,\alpha}}$.
\label{le:4.1}
\end{lemma}
The proof is straightforward. 

\subsection{CMC surfaces close to a half-Delaunay surface} 

We shall now study the problem of finding all CMC surfaces  near to a given 
half-Delaunay surface $D^+_\tau(s_0)$. These will be parametrized by their 
boundary values at $s=s_0$. As before,  we require some knowledge of the 
behavior of this solution as $\tau \to 0$. 

\subsubsection{The mean curvature operator on a Delaunay surface}

Using the isothermal parametrization (\ref{eq:2.1}), the unit normal  vector 
field $N_\tau$ on $D_\tau$ is given by 
\begin{equation}
N_\tau(s,\theta) :=( - \tau\,\cosh \sigma(s) \, \cos\theta, - 
\tau \, \cosh \sigma(s) 
\,\sin\theta, \del_ s \sigma(s)),
\label{eq:5.1}
\end{equation}
when $\tau \in (0, 1]$ and is given by 
\begin{equation}
N_\tau(s,\theta) :=(  \tau\,\sinh \sigma(s) \, \cos\theta, 
\tau \, \sinh \sigma(s) 
\,\sin\theta,  - \del_ s \sigma(s)), 
\label{eq:5.1bis}
\end{equation}
when $\tau \in (- \infty, 0)$. Normal graphs over $D_\tau$ then admit 
parametrizations of the form
\[
X_w: (s,\theta) \longmapsto X_\tau(s,\theta) + w(s,\theta) \, N_\tau(s,\theta),
\]
where $w$ is any function which is suitably small. 

\medskip

The mean curvature of the normal graph of a function $w$ about $D_\tau$ is 
calculated by a fairly complicated  nonlinear elliptic expression of $w$ which 
we shall not write out in full. However, this surface has mean curvature equal 
to $1$ if and  only if $w$ is a solution of the corresponding equation, 
which we write simply as $\calM_{\tau}(w) = 0$. We need  to know a bit 
about the structure of this operator, which is proved  in \cite{MP}~:
\begin{proposition}
The equation $\calM_\tau(w) = 0$ can be written in the form
\begin{equation} 
L_\tau w = \tau\, e^\sigma\, Q_\tau\left(\frac{w}{\tau\,e^\sigma }\right) 
\label{eq:5.2} 
\end{equation} 
where $L_\tau$ is the (modified) Jacobi operator (\ref{eq:3.1}), and  where 
$Q_\tau $ is a nonlinear second order differential operator which satisfies  
\[ 
Q_\tau (0) = 0 \qquad D_w Q_\tau (0) =0 \qquad \mbox{and}  \qquad  
D^2_w Q_\tau (0) =0.   
\] 
The Taylor expansion of $Q_\tau$ (in $w/\tau e^\sigma$) has coefficients 
which are uniformly bounded in $s$, along with all their derivatives,  
independently of $\tau \in (\tau_*,0) \cup  (0,1]$. 
\label{pr:5.1} 
\end{proposition}

We shall specialize now and suppose that $s_0 = -s_\tau$. Then, from  the 
discussion in \S 2.1, we see that the surface $D_\tau^+(s_\tau)$ (which we 
now simply call $D_\tau^+$) is  nearly flat close to its boundary, and it 
will be more convenient  computationally to use a slightly different  
parametrization of nearby surfaces, replacing the unit normal $N_\tau$ 
by a small perturbation of it, $\bar{N}_\tau$, which is constant (and 
in fact vertical, downward pointing) for all $s$ close to $-s_\tau$. 

\medskip

To define $\bar{N}_\tau$, choose a smooth cutoff function $\chi_\tau $ 
which is nonnegative and which satisfies $\chi_\tau  = 1$ when 
$s \geq -  s_\tau + 2$ and $\chi_\tau  = 0$ when $s \leq - s_\tau +1$. 
Now set
\[
\bar N_\tau : =  \chi_\tau \, N_\tau -  (1-\chi_\tau) \, (0,0,1).
\]
This satisfies
\begin{equation}
\left|\nabla^k \left(\bar N_\tau \cdot N_\tau - 1\right)   
\right|\leq c_k \, \tau, 
\qquad \mbox{for all}\quad k \geq 0 \quad \mbox{and}\quad s \in [-s_\tau, -
s_\tau+2]. 
\label{eq:6.1}
\end{equation}

We henceforth use the parametrization
\[
\bar{X}_w   : (s, \theta) \longmapsto X_\tau (s, \theta)  + w (s, \theta)  \, 
\bar{N}_\tau (s, \theta) , \qquad (s,\theta)\in [-s_\tau,\infty)  \times S^1.
\]
Denote by $D^+_\tau(w)$ the surface obtained this way.

\medskip

It follows from (\ref{eq:5.2}) and from (\ref{eq:6.1}) that $D^+_\tau(w
)$ has constant mean curvature equal to $1$ if and only if $w$ satisfies a  
nonlinear equation of the form
\[
L_\tau w = \bar{Q}_\tau (w),
\]
where 
\[
\bar{Q}_\tau(w)  := \ds \tau \,\bar L_\tau w + \tau\, e^{\sigma} \, 
\bar{Q}_{\tau}\,\left(\frac{w}{\tau\,e^{\sigma}}\right). 
\]
Here $\bar{Q}_{\tau}$ satisfies the same properties as listed for $Q_\tau$ in 
Proposition~\ref{pr:5.1}, and in fact these operators agree when $ s\geq - 
s_\tau+2$. The linear operator $ \tau \, \bar L_\tau$ represents  the deviation 
between the linearizations corresponding to the parametrizations $X_w$ using 
$N_\tau$ and  $\bar{X}_w$ using $\bar{N}_\tau$.  From (\ref{eq:6.1}), 
$\bar L_\tau$ has coefficients supported in  $[- s_\tau,- s_\tau+2] \times S^1$ 
which are  bounded in any  $\calC^{k,\alpha}$, uniformly in $\tau$. The details 
of this change of parametrization  are contained in \cite{MPP}.

\subsubsection{The nonlinear Poisson problem}

We are now ready to solve the nonlinear boundary problem $\calM_\tau(w) = 
0$  on $D_\tau^+$ with the value of $w$ at $s = -s_\tau$ (almost) specified.

\medskip

Fix $\mu \in (-2, -1)$ and $\delta >0$. Let $h$ be any element of  
$\calC^{2,\alpha}(S^1)$ which is orthogonal to $1$ and $e^{\pm i \theta}$, in 
the $L^2$ sense, and which satisfies
\[
\|h\|_{\calC^{2,\alpha}}\leq \delta \, \tau^{3}.
\]
Now define the approximate solution
\[
w_h : = \calP(h)(s + s_\tau ),
\]
where $\calP$ is the Poisson operator for $\Delta_0$ from 
Lemma~\ref{le:4.1}.  It suffices to use this Poisson operator, rather than the 
one for $L_\tau$,  because $L_\tau$ is very close to $\Delta_0$ in a long 
interval around  the boundary. Since we are using norms with exponential 
weight  factors,  these operators differ by a very small amount in norm. From 
the bounds in this lemma, the shift by $s_\tau$ in the $s$-variable,  and the 
fact $s_\tau \sim - \frac{1}{4}\log \tau^2$, we have 
\begin{equation}
\|w_h\|_{\calE^{2,\alpha}_{-2}} \leq c\, \tau \,\|h\|_{\calC^{2,\alpha}}.
\label{eq:6.2}
\end{equation}

We shall now search for a solution $w$, which we write as $w = w_h + v$, of 
$\calM_\tau(w) = 0$. The function $v$ will lie in  $\calE^{2,\alpha}_{\mu}([- 
s_\tau,\infty)\times S^1)$ and should satisfy
\[
\left\{ \begin{array}{rlll}
L_\tau v &  = &  \ds \bar{Q}_\tau (w_h + v)  - L_\tau w_h  \qquad & 
\mbox{in} \qquad (-s_\tau , \infty)\times S^1 \\[3mm]
v(-s_\tau,\cdot) &  \in   & \mbox{Span} \{ 1, e^{\pm i\theta} \}.
\end{array}
\right.
\]
The reason we are only requiring those eigencomponents of $v(s,\theta) 
=  \sum v_j(s)e^{ij\theta}$ with $|j| \geq 2$ to vanish 
is that we shall be using the 
inverse $G_{\tau,-s_\tau}$ from Proposition~\ref{pr:4.3},  which does not 
allow these low eigencomponents to be specified.  

\medskip

We solve this equation by finding a fixed point of the mapping
\[
\calN_\tau(v) :=  G_{\tau, -s_\tau }\, \left(\bar{Q}_\tau(w_h + v) - L_\tau 
w_h \right),
\]
>From (\ref{eq:6.2}) we have
\[
\|L_\tau w_h \|_{\calE^{0,\alpha}_{\mu}}\leq c\,\tau \, \|h 
\|_{\calC^{2,\alpha}}, \qquad \qquad   \|\tau \, \bar L_\tau w_h 
\|_{\calE^{0,\alpha}_\mu}\leq c \,  \tau^{1-\mu/2}\,\|h\|_{\calC^{2,\alpha}},
\]
and
\[
\|\tau\, e^\sigma\,\bar{Q}_{\tau}\,\left(\frac{w_h}{\tau\,e^\sigma} 
\right)\|_{\calE^{0,\alpha}_\mu} \leq c \,\tau^{-(3+\mu)/2}\,  \|h 
\|_{\calC^{2,\alpha}}^2.
\]
Only the final estimate uses that the norm of $h$ is small, and in fact,  
it is only 
really necessary to assume that $\|h\|_{\calC^{2,\alpha}} \leq c_0\,
\tau^{3/2}$ for $c_0>0$ sufficiently small.

\medskip

At this point it is straightforward, using Proposition~\ref{pr:4.4}, to 
show that 
there exists $c_* >0$ and $\tau_0 >0$, such that when  $|\tau| \in 
(0, \tau_0)$, 
the nonlinear mapping $\calN_\tau$ is a  contraction in the ball
\[
{\mathcal B} := \left\{v \quad | \quad \|v\|_{\calE^{2,\alpha}_\mu} \leq  c_* \, \tau  
\,\|h\|_{\calC^{2,\alpha}} \right\},
\]
and hence $\calN_\tau$ has a unique fixed point in this ball. Observe 
that $\tau_0 >0$ depends on $\delta$  while $c_* >0$ does not depend 
on $\delta$.

\medskip

In summary, we have proved
\begin{theorem}
Fix $\mu \in (-2, -1)$ and $\delta >0$. There are numbers $\tau_0 >0$ and 
$c>0$ such that for each $\tau$ with $0<|\tau| < \tau_0$ and  for every $h \in 
\calC^{2,\alpha}(S^1)$ which is orthogonal to $1$ and  $e^{\pm i \theta}$ and 
which satisfies $\|h\|_{2,\alpha}\leq \delta  \,\tau^3$, there exists an 
embedded CMC  surface $D_\tau(h)$, parameterized by 
\[
\bar{X}_w:= X_\tau + w \, {\bar N}_\tau\qquad \mbox{in}  \qquad  [- s_\tau, 
\infty ) \times S^1.
\]
The function $w$ here lies in $\calE^{2,\alpha}_\mu([-s_\tau,\infty)\times 
S^1)$ and satisfies $\|w\|_{\calE^{2,\alpha}_\mu} \leq c \, \tau^{-\mu/2} \,\| h 
\|_{\calC^{2,\alpha}}$ and finally, $w(-s_\tau,\cdot) - h(\cdot)  \in \mbox{\rm 
Span}\, \{1, e^{\pm i\theta}\}$.  
\label{th:6.1}
\end{theorem}

One of the main reasons for modifying the unit normal vector field  $N_\tau$ 
to $\bar{N}_\tau$ is because with this definition, the region near the boundary 
of $D_\tau(h)$ is a vertical graph over an annulus  with outer boundary $\del 
B_{r_\tau}$ in $\R^2$. Here $r_\tau = \frac{\tau}{2} e^{\sigma (- s_\tau)} \sim 
\tau^{3/2}$ is defined in \S 2.2. 

\medskip

We conclude this subsection with some estimates on the graph function in this 
representation. To do this, we first define some function spaces~:  
\begin{definition} 
For $r \in \N$, $\alpha \in (0,1)$ and $\mu \in \R$, we let  
$\calC^{r,\alpha}_{\mu}({\mathbb R} - 
\{0\})$ be the space of functions 
$w\in \calC^{r,\alpha}_{\mathrm loc}( {\mathbb R} - \{0\})$ such that 
\[
\|w\|_{\calC^{r,\alpha}_\mu} : = \sup_{\rho >0} \rho^{-\mu} \,  \|w ( 
\rho \, \cdot) \|_{\calC^{r,\alpha}(\overline{B_2}- B_1)} < \infty.
\] 
If $\Omega$ is a closed subset of ${\R}^2-\{0\}$, 
we define the space ${\mathcal C}^{r, \alpha}_{\mu}(\Omega)$ as 
the space of restriction of functions of  ${\mathcal C}^{r, \alpha}_{\mu} 
({\R}^2-\{0\})$ to $\Omega$. 
This space is naturally endowed with the induced norm. 
\label{de:6.1}
\end{definition}
Notice that the space $\calC^{r ,\alpha}_{\mu}(\overline{B_R} - \{0\})$ 
corresponds precisely to the space $\calE^{r ,\alpha}_\mu ([s_0,\infty)\times 
S^1)$ where $s_0 = -\log R$, under the  elementary change of variables $s 
= -\log r$.  However, we shall not use this obvious change of variables but 
rather the  more complicated change of variables specific to our problem  
\[
r := \frac{\tau }{2} \, e^{\sigma (s)},
\]
for $s \in [-s_\tau, 0)$ and $\theta \in S^1$.  Using the estimates of \S 2.2, 
we obtain 
\begin{equation}
r \, \del_r = ( 1 + {\cal O} (\cosh^{-2} s)) \, \del_s ,
\label{eq:cdv}
\end{equation}
for all  $s \in [-s_\tau, 0)$.

\medskip

We finally translate the surface $D_\tau (h)$ by $\mp  
\frac{\tau^2}{4} \, \log \left( \frac{4}{\tau^2}\right)$ along the 
vertical axis. This surface will still be  denoted by $D_\tau (h)$.
Using the expansion (\ref{eq:2.4}) together with (\ref{eq:cdv}), 
we see that  near $\del B_{r_\tau}$, the surface $D_\tau (h)$ 
is the graph of the function 
\begin{equation} 
\overline{B_{r_\tau}} - B_{r_\tau/2} \ni x \longmapsto \mp  
\frac{\tau^2}{4} \, \log \, r  - W_{h}(x) + V_{\tau, h}(x) ,
\label{eq:6.3}
\end{equation}
according to the sign of $\tau$ (with $-$ when $\tau >0$ and $+$ with $\tau 
<0$).  Here $W_{h}$ denotes the unique harmonic extension of $h$ in  
the ball $B_{r_\tau}$ and $V_{\tau, h} $ is bounded in $\calC_0^{2,\alpha} 
(\overline{B_{r_\tau}} - B_{r_\tau/2})$, which does not depend on $\delta$ 
nor on $\tau$, times $\tau^3$.  This constant $\tau^3$ has its origin in 
(\ref{eq:2.4}). 

\medskip

Observe that, reducing $\tau_0$ if this is necessary, we can assume that the 
mapping $h \rightarrow V_{\tau, h}$ is continuous and in fact smooth. With 
little work we also find that
\begin{equation}
\|V_{\tau , \tilde h} -V_{\tau, h}\|_{{\mathcal C}^{2, \alpha}_0}\leq c \,  
\tau^{1+\mu/2}  \, \|\tilde h - h \|_{{\mathcal C}^{2, \alpha}}
\label{eq:6.4}
\end{equation}
for some constant $c >0$ which does not depend on $\delta$, nor on $\tau$.

\section{The geometry and analysis of $k$-ended CMC surfaces}

\subsection{Moduli space theory for $k$-ended CMC surfaces}

We now briefly sketch the moduli space theory for $k$-unduloids as 
developed in \cite{KMP}. Because it is no harder to do so, we extend this 
and consider the deformation theory for complete, finite topology CMC 
surfaces with $k$ ends, each one of which is asymptotic to a Delaunay 
unduloid or nodoid with Delaunay parameter $\tau > \tau_*$, where 
$\tau_*$ is defined in (\ref{eq:deftst}). The statements and 
proofs are identical in this slightly broader context. 

\begin{definition} The moduli space $\calM_{g,k}^{\tau_*}$ consists
of the set of all complete constant mean curvature surfaces of finite
topology, with genus $g$ and $k$ ends $E_1, \ldots , E_k$ such that
each end $E_j$ is asymptotic to a half Delaunay surface
$D_\tau^+$ with $\tau \in (\tau_*,0)\cup(0,1]$.
\end{definition}

Now decompose each $\Sigma \in \calM_{g,k}^{\tau_*}$ into a union of
a compact component $K$ and ends $E_\ell$, $\ell=1, \ldots, k$. For each 
$\ell$ choose standard isothermal coordinates $(s,\theta)$ for the
model Delaunay end $D_{\tau_\ell}$, so that $E_\ell$ is parametrized by
\[ 
Y_\ell : =  X_\ell + w_\ell \, N_{\ell} \quad | \quad  [0, \infty)\times S^1
\longrightarrow E_\ell.
\] 
Here $X_\ell$ is the standard parametrization (\ref{eq:2.1}) 
for $D_{\tau_\ell}$;
each function $w_\ell$ decays exponentially and in fact
\begin{equation}
w_\ell \in \calE^{2,\alpha}_{-\gamma_{\tau_\ell,2}}([0,\infty)\times 
S^1). 
\label{eq:6.22}
\end{equation}
\begin{definition} 
For $r \in \N$, $\alpha \in (0,1)$ and $\mu \in \R$, let 
${\mathcal D}^{r,\alpha}_\mu(\Sig)$ be the space of functions 
$v \in{\calC}^{r,\alpha}(\Sig)$ for which 
\[ 
\| v \|_{\calD^{r,\alpha}_\mu} : = 
\|\left. v \right|_K\|_{\calC^{r,\alpha}} + 
\sum_{\ell=1}^k \|\left. v \circ Y_\ell \right|_{E_\ell}
\|_{\calE^{r,\alpha}_\mu }  < \infty. 
\] 
\label{de:7.1}
\end{definition}  

Let $\calL_\Sig = \Delta_\Sig + |A_\Sig|^2$ denote the Jacobi operator 
$\Sig$. Because of the asymptotic structure of the ends of $\Sig$,
the various mapping and regularity properties of this operator 
may be deduced from the analogous properties for $\calL_{\tau_\ell}$.
In particular, the set of indicial roots for $\calL_{\Sig}$, given by 
\[
\Gamma_\Sig := \{ \pm \, \gamma_{\tau_\ell, j} \quad | \quad j\in {\mathbb N}, 
\quad \ell=1, \ldots, k\},
\]
determine the weighted spaces on which $L_\Sig$ is Fredholm as well
as the asymptotic behavior of solutions of the homogeneous 
equation $\calL_\Sig w=0$. In particular, from the analysis of \S 2.5, 
when  $\mu \notin \Gamma_\Sig$, then
\[
\calL_\Sig : {\mathcal D}^{2,\alpha}_\mu(\Sig) \longrightarrow 
{\mathcal D}^{0,\alpha}_{\mu}(\Sig)
\]
is Fredholm. To keep track of the value of the weight parameter, we write 
this mapping as $\calL_\Sig(\mu)$.  As before, the operators 
$\calL_\Sig(\mu)$ and $\calL_\Sig(-\mu)$ are essentially dual to one 
another which implies that when $\mu \notin \Gamma_\Sig$, the operator 
$\calL_\Sig(\mu)$ is surjective if and only if the operator $\calL_\Sig(-\mu)$
is injective, and moreover $\mbox{dim\, ker}(\calL_\Sig(\mu)) = 
\mbox{dim\, coker}(\calL_\Sig (-\mu))$. 
We now give the precise definition of nondegeneracy.

\begin{definition} 
The surface $\Sig \in \calM_{g,k}^{\tau_*}$ is nondegenerate  if
$\calL_{\Sig}(\mu)$ is surjective for all $\mu > 0$ with 
$\mu \notin \Gamma_\Sig$, 
or equivalently, if $\calL_\Sig(-\mu)$ is injective for all $\mu>0$. 
\label{de:7.2}
\end{definition} 

\noindent Thus a surface is nondegenerate if it has no decaying Jacobi
fields. Note in particular that from the definition of $\tau_*$, the
Delaunay surfaces $D_\tau$ are nondegenerate when $\tau > \tau_*$. 

\medskip

The basic result from \cite{KMP} is
\begin{theorem}
Fix any element $\Sig \in \calM_{g,k}^{\tau_*}$. If $\Sig$ is nondegenerate,
then some neighborhood of $\Sig$ in the moduli space 
$\calM_{g,k}^{\tau_*}$ is a real analytic manifold of dimension $3k$. 
In general, the moduli space is a locally real analytic variety, 
i.e.\  there is a neighborhood of $\Sig \in \calM_{g,k}^{\tau_*}$
which is identified, via a real analytic diffeomorphism, with 
the zero set of a real analytic function defined 
in a finite dimensional Euclidean space.
\label{th:7.1}
\end{theorem}

We sketch the main ideas in the proof, 
but only in the nondegenerate case. The general case is somewhat more 
intricate, but is based on precisely the same ideas. 

The first point is that there exist analogues of the $6$ linearly 
independent geometric Jacobi fields $\Phi_{\tau_\ell}^{j,\pm}$,
$j = 0, \pm 1$, on each end $E_\ell$, $\ell = 1, \ldots, k$ of $\Sig$.
These are denoted $\Phi_\ell^{j,\pm}$.
\begin{lemma}
Let $\mu \in (-\gamma_{\tau_\ell,2},0)$.
Then these geometric Jacobi fields $\Phi^{j,\pm}_{\ell}$ satisfy
\[
\Phi^{j,\pm}_{\ell}\circ Y_\ell - \Phi^{j,\pm}_{\tau_\ell} 
\in \calE^{2,\alpha}_\mu(E_\ell).
\]
\end{lemma}
\noindent {\bf Proof :} Except for $\Phi_{\ell}^{0,-}$, these asymptotics 
may be deduced from the same constructions as in \S 2.4.1. The existence
and asymptotics of the remaining Jacobi fields $\Phi^{0,-}_{\ell}$ 
may be deduced by a perturbation argument using Proposition~\ref{pr:4.3} and 
(\ref{eq:6.22}). \hfill $\Box$

\medskip 

The proof of Theorem~\ref{th:7.1} is based on the implicit function theorem,
but to apply this theorem we need to find function spaces on which the
nonlinear mean curvature operator acts and on which $\calL_\Sig$ is 
surjective. Unfortunately, $\calL_\Sig(\mu)$ is never surjective
when $\mu < 0$. This may appear discouraging since the nonlinear operator 
as given only acts on spaces consisting of functions which decay along the
ends. To remedy this, we define the $6k$-dimensional deficiency space 
\[ 
{\mathcal W}_\Sig : =  \oplus_{\ell=1}^k \mbox{Span}\{\, \chi_\ell 
\,\Phi^{j, \pm}_{\ell} \quad | \quad j=-1, 0, 1\},
\] 
where $\chi_\ell$ is a cutoff function equal to $0$ near $E_\ell \cap
K$ and equal to $1$ on
$Y_{\ell}([c, \infty)\times S^1)$ for some $c>0$. 

\medskip

\begin{proposition} 
Assume that $\Sigma$ is nondegenerate and fix $\mu \in (-\inf_\ell 
\gamma_{\tau_\ell,2},0)$. Then the mapping 
\begin{equation}
\calL_\Sig :\calD^{2,\alpha}_\mu(\Sig) \oplus \calW_\Sig
\longrightarrow \calD^{0,\alpha}_\mu(\Sig)
\label{eq:7.1a}
\end{equation}
is surjective and its nullspace $\calN_\Sig$ is $3k$-dimensional.
Hence there exists a $3k$-dimensional subspace $\calK_\Sig \subset
\calW_\Sig$ such that
\begin{equation}
\calL_\Sig :\calD^{2,\alpha}_\mu(\Sig) \oplus \calK_\Sig
\longrightarrow \calD^{0,\alpha}_\mu(\Sig)
\label{eq:7.1}
\end{equation}
is an isomorphism.
\label{pr:7.1} 
\end{proposition}
The proof can be found in \cite{KMP}. 

\medskip

It is possible to make sense of the mean curvature operator on 
elements of the domain space in (\ref{eq:7.1a}) using the fact that 
elements of $\calW_\Sig$ correspond to geometric motions as in \S 2.4.1. 
Indeed, decomposing $u \in \calD^{2,\alpha}_\mu(\Sig)\oplus \calW_\Sig$ as
\[
u = u'+ u'',
\]
with $u'\in \calD^{2,\alpha}_\mu(\Sig)$ and $u''\in \calW_\Sig$
we then let $\Sig_u$ denote the surface which is the normal graph by
the function $u'$ over the surface obtained by slightly deforming 
the ends of $\Sig$ in the manner prescribed by the components of $u'' \in 
\calW_\Sig$. (More precisely, one defines a $6k$-dimensional real 
analytic parameter space $\calP$ of geometric deformations of $\Sig$ such 
that the differentials of curves in $\calP$ through $\Sig$ correspond
to the geometric Jacobi fields $\Phi_\ell^{j,\pm}$.) 
The mean curvature of $\Sig_u$ is identified with some function
$H(u)$ defined on $\Sig$, and by Proposition~\ref{pr:7.1}, when
$\Sig$ is nondegenerate, the differential of this map with respect to $u$ 
at $u=0$ is surjective. 

\subsection{CMC surfaces close to a $k$-unduloid}
 
\subsubsection{Geometric preparations}

Let $\Sig \in \calM_{g,k}^{\tau_*}$ be nondegenerate and fix a point 
$p_0 \in \Sig$. Assume by a rigid motion that $p_0 = 0$ and the 
oriented normal vector to $\Sig$ at that point is $(0,0,-1)$, so that 
the tangent plane $T_{0}\Sig$ is the $x \, y$-plane. Then in some 
neighborhood of $0$, $\Sig$ can be represented as a vertical graph
$z = u_0(x,y)$ for some function $u_0$ defined on a ball $B_\rho$.
By the assumptions above, $u_0(0) = \nabla u_0(0) = 0$. 

\medskip

Now for any sufficiently small vector $a = (a_{-1},a_0,a_1) \in \R^3$, 
we rotate and translate $\Sig$ slightly so that this new surface $\Sig_a$
is the graph of a function $u_a$, also defined on $B_\rho$, such that
\[
u_a(0) = a_0, \qquad \del_x u_a (0) =  a_1, \qquad 
\mbox{and} \qquad \del_y u_a (0) =  a_{-1}.
\]
Also, for $0 < r \leq  \rho$, let $\Sig_{a,r}$ denote the complement in 
$\Sig_a$ of the graph of $u_a$ on $B_r$. Finally, set 
\[
p_a : =  (0, 0,u_a (0)).
\]

\medskip

Next, let us modify the unit normal $N_{a}$ on $\Sig_a$ to a new
unit vector field $\bar N_a$ on $\Sig_a$ such that 
in $\Sig_a - \Sig_{a,\rho/2}$,  ${\bar N}_{a} \equiv
(0,0,-1)$, while on $\Sig_{a,\rho}$, $\bar N_a = N_\Sig$.
The linearization of the mean curvature operator with respect to this 
new vector field, $\bar{\calL}_a$, is a slight perturbation of 
$\calL_{\Sig_a}$. Any mapping property for $\calL_\Sig$ immediately
transforms to one for $\bar{\calL}_a$, hence in particular 
Proposition~\ref{pr:7.1} holds when $\calL_\Sig$ is replaces by 
$\bar{\calL}_a$.

\subsubsection{The mean curvature operator for graphs}

The mean curvature operator for the vertical graph of a function $u: 
\R^2 \supset B_\rho \rightarrow \R$ (with downward pointing normal)
is given by 
\[
H(u) :=  - \frac12 \mbox{div} \, \left(\frac{\nabla u}{(1+|\nabla 
u|^2)^{1/2}}\right).
\]
Since $\bar{\calL}_a \, w = \left. D_u H\right|_{u=u_a}(w)$ in
$\Sig_a - \Sig_{a,\rho/2}$, we obtain that 
\[
\bar{\calL}_a w  =  \frac12 \mbox{div} \, \left(\frac{\nabla w}{
(1+|\nabla u_a|^2)^{1/2}}\right)  - \frac{1}{2} \, \mbox{div} \left(   
\frac{\nabla w \cdot \nabla u_a }{(1+|\nabla u_a|^2)^{3/2}} 
\, \nabla u_a \right)
\]
there. The main point is that $\bar{\calL}_a$ is close to the 
standard Laplace operator. More precisely, we have the
\begin{lemma}
Fix any $\nu \in \R$. There exist $a_*,c >0$ such that when $|a_{\pm 1}| 
< a_*$, then for $0 < \rho_1 < \rho_2$ and any 
$w \in \calC^{2,\alpha}_\nu(\overline{B_{\rho_2}} - B_{\rho_1})$,
we have 
\[
\|\bar{\calL}_a w - \Delta w \|_{\calC^{0,\alpha}_{\nu -2} 
(\overline{B_{\rho_2}} - B_{\rho_1})}  \leq c \, 
(\rho_2 \,(|a_1| + |a_{-1}|) + \rho_2^2 ) \, 
\|w\|_{\calC^{2,\alpha}_\nu(\overline{B_{\rho_2}} - B_{\rho_1})}, 
\]
\label{le:klm}
\end{lemma}
A brief calculation also shows that the equation $H(u) = 1$ can be written as
\begin{equation}
\Delta u + \left(\Delta u \, |\nabla u|^2 - \frac12 \, \nabla^2 u \, 
(\nabla u,\nabla u)\right) + 2\, \left( (1 +|\nabla u|^2)^{3/2}-1-\frac32 
\, |\nabla u|^2 \right) =1.
\label{eq:mlml}
\end{equation}

\subsubsection{Mapping properties of $\bar{\calL}_a$}

The proof of the following result can be found in \cite{FP}~:
\begin{lemma}
Fix $-1 < \nu < 0$ and $0 < \rho_1 < \rho_2/2$. Then there exists an operator 
\[
G_{\rho_1,\rho_2}: \calC^{0,\alpha}_{\nu-2}(\overline{B_{\rho_2}}
-  B_{\rho_1}) \longrightarrow \calC^{2,\alpha}_{\nu}
(\overline{B_{\rho_2}}- B_{\rho_1})
\]
with the following properties. For any $f \in \calC^{0,\alpha}_{\nu - 2}
(\overline{B_{\rho_2}} - B_{\rho_1})$, the function $w= G_{\rho_1, 
\rho_2}(f)$ is a solution of the problem 
\[
\Delta w =  f \qquad  \mbox{in} \quad B_{\rho_2}-
\overline{B_{\rho_1}} 
\]
with boundary data $w =0$ on $\del B_{\rho_2}$ and $w$ equal to a constant
on $\del B_{\rho_1}$. In addition, 
\[
\|G_{\rho_1, \rho_2}(f) \|_{\calC^{2,\alpha}_\nu} 
\leq c \, \|f \|_{\calC^{0,\alpha}_{\nu - 2}}
\]
for some constant $c>0$ independent of $\rho_1$ and $\rho_2$.
\label{le:8.01}
\end{lemma}

We now define the function spaces suitable for the analysis of these operators
on the punctured surfaces $\Sig_a  - setminus \{p_a\}$. We identify all
functions defined near $p_a$ on $\Sig_a$ with functions on $B_\rho \subset
\R^2$ using the fixed coordinates induced by the functions $u_a$. 
\begin{definition} 
Let $r \in \N$, $0 < \alpha<1$ and $\nu, \mu \in \R$. Then 
$\calD^{r,\alpha}_{\mu,\nu}(\Sig_{a}  -  \{p_a\})$ is 
the space of functions $w \in \calC^{r,\alpha}_{{\mathrm{loc}}}
(\Sig_a  -  \{p_a\})$ for which 
\[
\|w\|_{\calD^{r,\alpha}_{\mu,\nu}} : =  \|w\|_{\calD^{r,\alpha}_\mu 
(\Sig_{a,\rho})} + \| w \circ u_a^{-1}\|_{\calC^{r,\alpha}_\nu 
(B_\rho  -  \{0\})} < \infty.
\] 
\label{de:8.1}
\end{definition}
Thus $\mu$ is the weight on the other ends of $\Sig_a$ and $\nu$
regulates growth or decay near $p_a$. 
As usual, if $\Omega$ is any closed subset of $\Sig_a  -  \{p_a\}$,
then $\calD^{r,\alpha}_{\mu,\nu}(\Omega)$ is the space of 
restriction of functions in $\calD^{r,\alpha}_{\mu,\nu}(\Sig_a  - 
\{p_a \})$ to $\Omega$, endowed with the induced norm. 

\medskip

\begin{proposition}
Fix $\mu \in (- \inf_\ell \gamma_{\tau_\ell,2}, 0)$, $-1 < \nu <0$ and 
$0 < \alpha <1$. Recall the number $r_\tau$ defined in $\S\ref{singlim}$. 
Then for $\tau$ in some small interval $(0,\tau_0)$, 
there exists an operator 
\[
\bar{\mathcal G}_{\tau}:\calD^{0,\alpha}_{\mu,\nu-2} 
(\Sig_{a,r_\tau}) \longrightarrow \calD^{2,\alpha}_{\mu,\nu} 
(\Sig_{a,r_\tau}) \oplus {\mathcal K}_\Sig,
\]
such that for all $f \in \calD^{0,\alpha}_{\mu,\nu -2}(\Sig_{a,r_\tau})$, 
the function $w = \bar{\mathcal G}_{\tau}(f)$ is a solution of 
\[
\bar{\calL}_a w = f
\]
in $\Sig_{a,r_\tau}$, with $w$ constant on $\del \Sig_{a,r_\tau}$. 
Furthermore, there exists $c>0$ which is independent of $\tau$ and $f$,
such that
\[
\|w\|_{\calD^{2,\alpha}_{\mu,\nu}\oplus {\mathcal K}} \leq c \, 
\|f\|_{\calD^{0,\alpha}_{\mu,\nu-2}}. 
\]
\label{pr:8.1}
\end{proposition}
\noindent {\bf Proof:} First choose $\rho_* \in (0,\rho]$ and define
\[
w := \chi \, G_{\rho_*,r_\tau}(f),
\]
where $G_{\rho_*,r_\tau}$ is the operator in Lemma~\ref{le:8.01} and where 
$\chi$ is a radial cutoff function identically equal to $1$ in 
$B_{\rho_*/2}$ and equal to $0$ outside $B_{\rho_*}$. By construction the
support of the function $g$ defined by $g : = \Delta w -f$ 
is disjoint from $B_{\rho_*}$. Now, using Proposition~\ref{pr:7.1}, we 
define 
\[
w' := \chi' \, \bar{\calL}_a^{-1}(g),
\]
where $\chi'$ is another radial cutoff function equal to $1$ outside 
$B_{2r_\tau}$ and vanishing in $B_{r_\tau}$. 

\medskip

It is easy to check that
\[
\|f - \bar{\calL}_a \, w'\|_{\calD^{0,\alpha}_{\mu,\nu-2}} \leq 
c \, (\rho_*^2  + |a_{\pm 1}|\, \rho_* + r_\tau^{-\nu} ) \, \|f\|_{\calD^{0,
\alpha}_{\mu,\nu-2}},
\]
and also that 
\[
\|w'\|_{\calD^{2,\alpha}_{\mu,\nu}}\leq c \, \|f\|_{\calC^{0,\alpha}_{\mu,
\nu-2}} ,
\]
for some constant $c>0$ which depends neither on $\tau$ nor on $\rho_*$. 
The result follows immediately by a simple perturbation argument, 
provided $\rho_*$ and $\tau$ are chosen small enough. \hfill $\Box$

\medskip

We conclude this subsection with a simple result regarding the Poisson
operator for $\Delta$ which is close to 
Lemma~\ref{le:4.1} in spirit. 

\begin{lemma} 
For any $g \in \calC^{2,\alpha}(S^1)$ such that $\ds \int_{S^1} g = 0$,
there exists a unique function $w \in \calC^{2,\alpha}_{-1}({\R}^2
 -  B_1)$ which is a solution of
\begin{equation}
\left\{ 
\begin{array}{rllll} 
\Delta w & =  & 0 \qquad & \mbox{\rm in}\quad {\R}^2  - 
B_1 \\[3mm] 
w & =  &  g \qquad & \mbox{\rm on} \quad \del B_1,
\end{array}  
\right. 
\label{eq:8.0}
\end{equation}
and which satisfies $\| w \|_{\calC^{2,\alpha}_{ -1}} \leq c \, \| g 
\|_{\calC^{2,\alpha}}$ for some constant $c>0$ which does not depend on $g$. 
\label{le:8.1}
\end{lemma}

We write this solution of (\ref{eq:8.0}) as $P(g)$. 

\subsubsection{The nonlinear Poisson problem}

Using Proposition~\ref{pr:7.1} let $\gamma_0$ be the solution of
\begin{equation}
\bar{\calL}_a \gamma = -  2 \, \pi \, \delta_0, \qquad \mbox{in} \qquad 
\Sig_a  -  \{ p_a \},
\label{eq:8.1}
\end{equation}
with $\gamma + \chi \, \log r  \in \calD^{2,\alpha}_{\mu}(\Sig_a) 
\oplus {\mathcal K}_{\Sig_a}$, where $\chi$ is, as usual, a cutoff 
function equal to $1$ in $\Sig_a  -  \Sig_{a,\rho/2}$ and 
vanishing in $\Sig_{a, \rho}$ and $r=|(x,y)|$. 

\medskip

We now observe that there are three global Jacobi fields $\Phi_j$,
$j = 0, \pm 1$, on $\Sig_a$ such that
\[
\Phi_0 = 1 + {\mathcal O} (r), \qquad
\Phi_{1} = r \, \cos \theta + {\mathcal O} (r^2),\qquad 
\mbox{and}\qquad 
\Phi_{-1} = r \, \sin \theta + {\mathcal O} (r^2)
\]
in $B_\rho$. Indeed, $\Phi_0$ is obtained by projecting the constant 
Killing field $(0,0,1)$ (corresponding to a vertical translation)
onto the normal vector field, while $\Phi_{\pm 1}$
are obtained by projecting the Killing fields $(z,0,-x)$ and $(0,z,-y)$
(corresponding to the two rotations about the horizontal axes).
By adding an appropriate linear combination of these three Jacobi
fields to $\gamma$, we may also assume that $\gamma + \chi \log r$ 
and its gradient vanishes at $0$. 
\begin{lemma}
If $\gamma$ is defined as above, then for all $k \geq 0$, there exists a
constant $c_k > 0$ such that
\begin{equation}
\left| \nabla^k \left( \gamma + \log r \right) \right| \leq  c_k \, r^{2-k} \, |\log r|
\label{eq:8.2}
\end{equation}
in $B_{\rho}$, where $r := |x|$.
\label{le:8.2}
\end{lemma}
The proof is straightforward. 

Now fix $\mu \in (-\inf_\ell \gamma_{\tau_\ell,2}, 0)$, $-2/3 < \nu < 0$ and 
$\delta >0$. Suppose that $a_j \in \R$, $j=0,\pm 1$, satisfy
\[
|a_0|+ \tau^{3/2} \, |a_{\pm 1}|\leq \delta \, \tau^3.
\]
Finally, for any $g \in \calC^{2,\alpha}(S_1)$ such that $\int_{S^1} g = 0$,
and $\|g\|_{\calC^{2,\alpha}}\leq \delta \, \tau^3$, we define 
\[
w_g := P(g)(r_\tau^{-1}u_a^{-1}(\cdot)).
\]
This is the bounded harmonic extension of $g$ to $\R^2  -  B_{r_\tau}$.
By Lemma~\ref{le:8.1} there exists a constant $c>0$ such that 
\begin{equation}
\|w_g\|_{\calC^{2,\alpha}_{\mu,-1}}\leq c \, r_\tau \,  
\|g\|_{\calC^{2,\alpha}}.
\label{eq:8.22}
\end{equation}
Now define
\[
\tilde{w}_{\tau, g}  := \chi \, w_g - \frac{\tau^2}{4} \, \gamma ,
\]
where $\chi$ is the same cutoff function used in the definition
of $\gamma$. 

\medskip

We wish to find $v \in \calD^{2,\alpha }_{\mu,\nu}(\Sig_{a,r_\tau})
\oplus {\mathcal K}_{\Sig_a}$ so that the graph over $\Sig_a$ of the 
function $w := \tilde{w}_{\tau,g} + v$, using the 
vector field $\bar{N}_a$ is CMC. (Note that, just as in \S 3.1, 
we do not actually consider the graph of $w$ over $\Sig_{a, r_\tau}$,
but rather decompose $w = w'+ w''$ and consider the graph of $w'$ over 
some deformation $\Sig_{a,r_\tau,w''}$ of $\Sig_{a,r_\tau}$.)
This is equivalent to finding the solution of some nonlinear elliptic 
operator which we write formally simply as
\[
\bar{\mathcal L}_a v  =  Q(\tilde w_{g, \tau} +  v) - \bar{\calL}_a 
\tilde w_{g,\tau}.
\]
As usual, we use Proposition~\ref{pr:8.1} to rephrase this as a fixed 
point problem for the operator 
\[
{\mathcal M}_\tau(v) : = \bar{\mathcal G}_{\tau}\, (Q(\tilde w_{g, 
\tau} + v)-\bar{\calL}_a \tilde w_{g,\tau}).
\]

>From the estimates (\ref{eq:8.2}) and (\ref{eq:8.22}), along with
(\ref{eq:mlml}), it is easy to show that 
\[
\|\bar{\calL}_a \tilde w_{g, \tau}\|_{\calD^{2,\alpha}_{\mu,\nu - 2}}
\leq c \, \tau^{4}, \qquad \mbox{and} \qquad 
\|Q(\tilde w_{g,\tau})\|_{\calD^{2,\alpha}_{\mu,\nu - 2}}\leq   
c_\delta \,\tau^{9/2},
\]
where $c>0$ does not depend on $\delta$, provided $\tau$ is small enough. 

\medskip

It is then routine to show that for any fixed $\delta >0$, there exist 
$c_* >0$ and $\tau_0 >0$ such that when $\tau \in (0, \tau_0)$, the 
nonlinear mapping ${\mathcal M}_\tau$ is a contraction mapping in the 
ball
\[
\widehat{\mathcal B} :=  \{v \quad | \quad ||v||_{\calD^{2,\alpha}_{\mu,\nu} 
\oplus {\mathcal K}}\leq c_* \, \tau^{4} \},
\]
and thus has a unique fixed point in this ball. Again, $\tau_0 $ depends 
on $\delta$ while $c_*$ does not.

\medskip

In summary, we have proved the
\begin{theorem}
Fix constants $\mu \in (-\inf_\ell \gamma_{\tau_\ell,2}, 0)$, 
$\nu \in (-2/3, 0)$, $\alpha \in (0, 1)$ and $\delta >0$. Suppose that
$\tau$ sufficiently small and $a \in \R^3$ with $|a_0| + \tau^{3/2}\, 
|a_{\pm 1}|\leq \delta \, \tau^{3/2}$. Then for any $g \in \calC^{2,
\alpha}(S^1)$ with $\int_{S^1}g = 0$ and $\|g\|_{\calC^{2,\alpha}}\leq 
\delta \, \tau^3$, there exists a CMC surface with boundary which is 
close to $\Sig_{a,r_\tau}$ and such that a collar neighbourhood of its
boundary can be parameterized as a vertical graph
\begin{equation}
\overline{B_{2r_\tau}} -  B_{r_\tau} \ni x  \longrightarrow  - 
\frac{\tau^2}{4} \log r + a_0 +  a_1\, r \,  \cos\theta + a_{-1}\sin\theta
- \widehat W_g (x) + \widehat V_{\tau,a,g}(x).
\label{eq:8.3}
\end{equation}
Here $\widehat W_g$ is the unique bounded harmonic extension of $g$ 
to $\R^2  -  B_{r_\tau}$ and $\widehat V_{\tau,a,g}$ is a function 
which is bounded in $\calC_0^{2,\alpha}(\overline{B_{2r_\tau}} - 
B_{r_\tau})$ by $c \tau^3$, where $c>0$ does not depend on $\delta$ or 
$\tau$.
\end{theorem}
\noindent It is only in the final estimate that we need to restrict
$\nu$ to lie in $(-2/3, 0)$.  

\medskip

We can also see from this that $\widehat V_{\tau, a, g}$ depends
smoothly on the parameters $(\tau,a,g)$. In fact, with a little
work, we find that
\begin{equation}
\|\widehat V_{\tilde \tau,\tilde a , \tilde g }-V_{\tau,a,g}\|_{\calC^{2,
\alpha}_0}\leq c\,\tau^{3 \,\nu/2} \, ( \tau^{3/2}\, \|\tilde g-g 
\|_{\calC^{2,\alpha}}+\tau^4 \, |( \tilde{a}_{-1}, \tilde{a}_1 ) - 
(a_{-1},a_1)| + \tau^{3} \, |\tilde\tau -\tau|)
\label{eq:8.4}
\end{equation}
for some constant $c >0$ which not depending on $\delta$ or $\tau$.

\begin{remark}
A similar construction holds {\it mutatis mutandis} if one replaces 
$\gamma$ by $- \gamma$ in the formula defining ${\tilde w}_{\tau, g}$.
We denote by $\Sig^\pm(\tau, a, g)$ the surface constructed here,
where for $\Sig^+$ we have used $+\gamma$ and for $\Sig^-$ we have used
$-\gamma$ in the construction.
\end{remark}

\section{Adding a Delaunay end}

We now assemble the results established in the previous sections
and prove the main gluing theorem

\begin{theorem} Let $\Sig \in \calM_{g,k}^{\tau_*}$ be nondegenerate.
Then there exists $\tau_0 > 0$, ${\cal U}$ a neighborhood of $p$ in 
$\Sigma$ and a smooth gluing map
\[
{\mathfrak G}_\Sigma : {\cal U} \times \{\tau \quad | \quad 0 < |\tau| < \tau_0\} \longrightarrow
\calM_{g,k+1}^{\tau_*}
\]
where $(p,\tau)$ is mapped to the surface obtained by gluing
$D_\tau^+$ to $\Sigma$ at $p$. The surfaces in the range of this
map ${\mathfrak G}$ are all nondegenerate.
\label{th:mgt}
\end{theorem}
We prove this theorem in the next two subsections.

\subsection{Step 1~: existence} 

Let $\Sig \in \calM_{g,k}^{\tau_*}$ be nondegenerate and fix a point $p \in 
\Sig$ and $\delta >0$ sufficiently large enough, and apply the results 
of the previous sections. These give the existence of a $\tau_0 > 0$ 
such that when $0 < |\tau| < \tau$, there are two families
of CMC surfaces with boundary~: $D_\tau^+(h)$ and $\Sig^\pm(\tilde{\tau},a,g)$.
The first of are the perturbed half-Delaunay surfaces, 
where the boundary values 
$h \in \calC^{2,\alpha}(S^1)$ satisfies 
\[
\int_{S^1} h = \int_{S^1} h e^{\pm i\theta} = 0, 
\qquad \mbox{and}\qquad
\|h\|_{2,\alpha} \leq \kappa \tau^3.
\]
The second are the CMC perturbations of $\Sig_{a,\rho/2}$, where the
parameters are defined as follows.  First,
\[
\tilde \tau : = \sqrt{\tau^2 + 4 \,  t},
\]
where $t\in \R$ is any number such that 
\[
|\log \tau| \, |t|\leq\delta \,\tau^3;
\]
next $a \in \R^3$ satisfies
\[
|a_0|+\tau^{3/2} \, |a_{\pm 1}| \leq \delta \, \tau^3;
\]
finally, $g \in \calC^{2,\alpha}(S^1)$ is such that
\[
\int_{S_1}g = 0 \qquad \mbox{and} \qquad 
\|g\|_{2,\alpha}\leq \kappa \,\tau^3.
\]
The superscript on $\Sig^\pm(\tilde{\tau},a,g)$ should be taken to be
$+$ when $\tau >0$ and $-$ when $\tau < 0$. On the other hand 
the $+$ superscript on $D_\tau^+(h)$ only serves to remind us 
that we are dealing with a half-Delaunay surface, since $\tau$ 
(and hence its sign) appears explicitly. 

\medskip

If we can now find choices of $g$, $h$, $t$ and $a$ so that 
\[
\Sig^\pm(\tilde\tau,a,g) \sqcup D_\tau^+(h) 
\]
is $\calC^1$ across the joining interface, then the existence will
be established. This is because the equation is noncharacteristic
at the boundaries where these surfaces intersect and standard regularity 
theory for the mean curvature equation shows that this union is then 
a $\calC^\infty$ surface. 
We denote the CMC surfaces obtained in this way by ${\mathfrak G}_\Sigma (p, \tau)$.
These look like $\Sig$ with an additional Delaunay end,  
unduloidal when $\tau >0$ and nodoidal when $\tau <0$, attached at
the point $p \in \Sig$. 

\medskip

We now show that we can match the Cauchy data by choosing these
parameters correctly. By construction, near each of their respective
boundaries, the surfaces $\Sig^\pm(\tilde\tau,a,g)$ and $D_\tau^+(h)$ 
are vertical graphs over the $x \,y$-plane. Therefore, in view of
(\ref{eq:6.3}) and (\ref{eq:8.3}), it suffices to solve the equations
\begin{equation}
\left\{
\begin{array}{rllll}
\ds \mp  t \log r + a_0 + a_1\, r \, \cos\theta + a_{-1}\,r \,\sin 
\theta  - \widehat W_{g} + \widehat V_{\tilde\tau,a,g} & = & 
\ds - W_{h}+ V_{\tau, h} \\[3mm]
\ds \mp  t + a_1 \, r \,  \cos\theta + a_{-1}\, r \, \sin\theta  -  
r\, \del_r\widehat W_{g} + r \, \del_r\widehat V_{\tilde\tau,a,g} & = & 
\ds - r\, \del_r W_{h}+ r\,\del_r V_{\tau, h}; 
\end{array}
\right.
\label{eq:9.1}
\end{equation}
all functions here are evaluated on $\del B_{r_\tau}$. These identities
correspond to the coincidence of the Dirichlet and Neumann data,
respectively, of the two graphs. 

\medskip

To solve these, recall that the mapping 
\[
{\cal P} :  \calC^{2,\alpha}(S^1) \ni h \longrightarrow r_\tau \, 
\del_r (W_h - \widehat W_h)(r_\tau \, \cdot) \in \calC^{1,\alpha} 
(S^1)
\]
is an isomorphism such that both it and its inverse have norm
uniformly bounded in $\tau$. To see this first observe that 
this mapping does not depend 
on $r_\tau$, so we may as well assume that $r_\tau \equiv 1$. Next, note that
${\cal P}$ is a linear first order elliptic self-adjoint pseudodifferential 
operator with principal symbol $ - 2 \, |\xi|$. Thus, it is enough to
check that it is injective. Now, if ${\cal P} (h)=0$ then the function 
$w$ which equals $\widehat W_h$ in $\R^2  -  B_1$ and $W_h$ in 
$B_1$ is a global solution of $\Delta w =0$ on all of $\R^2$; 
furthermore, $w$ belongs to $\calC^{2,\alpha} (B_1) \cap 
[\calC^{2,\alpha}_{-1}(\R^2  -  B_1)\oplus \mbox{Span}\{\log r \}]$. 
No such function exists, and so first $w$ and then $h$ must be trivial.

\medskip

Now define 
\[
h^- : = g  \pm  t \, \log r_\tau \qquad 
\mbox{and} \qquad
h^+ : = h + a_0 + a_1 \, r_\tau \, \cos\theta + a_{-1} \, r_\tau \,  
\sin \theta.
\]
It is easy to see that (\ref{eq:9.1}) reduces to a fixed point problem 
\[
(h^+, h^-) = {\bf C}_\tau (h^+, h^- )
\]
in the space $\calE: = (\calC^{2,\alpha}(S^1))^2$.  But 
(\ref{eq:6.3}) and (\ref{eq:8.4}) imply that ${\bf C}_\tau : \calE 
\longrightarrow \calE$ is a contraction mapping defined in the ball of 
radius $\delta \, \tau^3$ in $\calE$ into itself, provided $\delta$ is 
sufficiently large $\tau$ is sufficiently small. This gives the 
fixed point and completes the proof of the existence of the gluing map. 

\subsection{Step 2~: nondegeneracy}

We now prove that the surfaces ${\mathfrak G}_\Sigma (p, \tau)$ 
constructed above
are nondegenerate when $\tau$ is small enough. The proof is by 
contradiction. Assume that this is not the case, so that there
exists a sequence $\tau_n \to 0$, a sequence $p_n \in \Sigma$ 
tending to $p$ and a sequence of surfaces  
\[ 
\Sig_n := {\mathfrak G}_\Sigma (p_n, \tau_n) ,
\] 
for which the operators $\calL_{\Sig_n}$ are not injective on 
$\calD^{2,\alpha}_{\mu_n}(\Sig_n)$, for some $\mu_n <0$. In other words, 
for each $n$ there exists a nontrivial function $w_n \in 
\calD^{2,\alpha}_{\mu_n}(\Sig_n)$ such that $\calL_{\Sig_n} w_n= 0$. 
Without loss of generality, we can assume that $p=0$ and that the 
tangent plane of $\Sigma$ at $0$ is the $x\, y$-plane.

\medskip 

The set indicial roots of $\calL_{\Sig_n}$ decomposes into two groups~:
\begin{itemize}
\item[(i)] Those indicial roots associated to the ends of 
$\Sig_n$ which converge to the ends of $\Sig$;
\item[(ii)] Those indicial roots associated to the end of $\Sig_n$ 
which is a perturbation of $D_{\tau_n}$.
\end{itemize}
As $n \to \infty$, elements of the first subset converge to the 
corresponding indicial roots of $\Sig$, while elements of the 
second subset converge to elements in ${\mathbb Z}  -  \{\pm 1\}$. 
Hence there exist $\mu  \in (-\inf_\ell \gamma_{\tau_\ell,2}, 0)$
and $\tilde{\mu} \in (-2,-1)$ such that when $n$ is large, $w_n$ is 
bounded by a multiple of $e^{\mu s}$ on each end of the `original'
ends of $\Sig_n$ and by a multiple of $e^{\tilde{\mu}s}$ 
on the new end of $\Sig_n$. 

\medskip 

Choose $\mu' \in (\mu,0)$ and $\tilde{\mu}' \in (\tilde{\mu}, -1)$. 
By construction, $\Sig_n$ decomposes into the union of surfaces,
one of which we denote by $\tilde\Sig_n$ and has boundary 
$\partial\tilde\Sig_n$ which is a normal 
graph a compact portion of $\Sig$. The other, which we denote by
$\tilde D_n$, is a normal graph over the Delaunay surface $D_{\tau_n}^+$. 
Now define a weight function $\zeta_n >0$ on $\Sig_n$, such that
\begin{itemize} 
\item $\zeta_n \sim e^{\mu s}$ on each end of $\tilde \Sig_n$ , 
\item $\zeta_n \sim r^{- \tilde{\mu}}$ near $\del \tilde \Sig_n$ , 
\item $\zeta_n \sim \tau^{- 2\tilde{\mu}} \, e^{\tilde{\mu}s}$ on  
$\tilde D_n$. 
\end{itemize} 
Here $f \sim g$ means that $1/2 \leq f/g\leq 2$. We use the
usual cylindrical coordinates to parametrize the various ends of 
$\tilde \Sig_n$, while in a small annulus near $\del\tilde \Sig_n$ 
we use polar coordinates. 

\medskip 

Normalize the sequence $w_n$ so that  
\[ 
\sup_{\Sig_n} \, \zeta_n^{-1}\, w_n = 1. 
\] 
By the choices of $\mu$ and $\tilde{\mu}$, these suprema are achieved
at some point $q_n \in \Sig_n$. We distinguish a few cases according to 
the behavior of the sequence of points $q_n$. Note that as $n \to \infty$,
the surfaces $\Sig_n$ converge to the union of the original surface
$\Sig$ and an infinite union of spheres of radius $1$ centered at the 
points $(0,0,2 \, j + 1)$, for $j \in {\mathbb N}$. 

\medskip 

\noindent 
{\em Case 1~:} Suppose that (some subsequence of) the points $q_n$ 
remain in a fixed end of $\tilde\Sig_n$ and tend to infinity. Then 
$q_n = (s_n,\theta_n) \in [0, \infty)\times S^1$, where 
$s_n \to \infty$. Possibly extracting another subsequence, 
we may assume that 
\[ 
\tilde w_n  := e^{-\mu s_n}\, w_n (s_n + \cdot, \cdot)  
\] 
converges uniformly on any compact subset of $\R\times S^1$ to 
a limiting function $w_\infty$ which is a solution of  
\[ 
\calL_{\tau_\ell} w_\infty =0 
\] 
on the full Delaunay surface $D_{\tau_\ell}$. Here $\tau_\ell$ is the
Delaunay parameter for that end in the surface $\Sig$. In addition,
$|w_\infty| \leq c \, e^{\mu s}$. To see this is impossible, 
decompose $w_\infty = \sum_j w_j  \, e^{i \, j \, \theta}$. By the
choice of $\mu$ we have $w_0 = w_{\pm 1} =0$ since the nontrivial
solutions of $\calL_{\tau_\ell} w=0$ in these eigenspaces are either 
bounded or blow up linearly at both ends of $D_{\tau_\ell}$. 
On the other hand, the restriction of $\calL_{\tau_\ell}$ to 
the eigenspaces with $|j| \geq 2$ satisfies the maximum principle.
This implies  that $w_j = 0$ for all other values of $j$. This is
a contradiction.

\medskip 

\noindent 
{\em Case 2~:} Next, suppose that the sequence $q_n$ converges to a
point $q_\infty \in \Sig  -  \{0\}$.  We may clearly assume that 
$w_n$ converges uniformly on any compact of $\Sig  -  \{0\}$ to 
a solution of  
\[ 
\calL_\Sig w_\infty =0 
\] 
defined on $\Sig -  \{0\}$. But $|w_\infty| \leq c \, 
r^{- \tilde \mu}$ as $r \to 0$ and $|w_\infty|\leq c\, e^{\mu s}$ at 
all other ends of $\Sig$. Hence it is a Jacobi field which decays 
exponentially at all ends, and so must vanish by nondegeneracy. 

\medskip 

\noindent 
{\em Case 3~:} Finally, suppose that $p_n$ tends to a point on the union 
of spheres centered at the points $(0,0,2j+1)$, for $j \in {\mathbb N}$. 
But $p_n$ corresponds to $(s_n, \theta_n)$ in the parameterization of 
Delaunay end of parameter $\tau_n$. Both $\tilde D_n$ and 
$\tilde \Sig_n  -  \tilde \Sig_{n,\rho}$ are normal graphs over 
$D_{\tau_n}$ when $\rho$ is small enough. Thus we may define 
the rescaled sequence
\[ 
\tilde{w}_n (s, \theta) : = e^{\mu \, s_n}\, w_n ( s + s_n, \theta). 
\] 
It is proved in \cite{MP} that, as $\tau_n \to 0$, 
the term of order $0$ in $\calL_{\tau_n}$ converges either to $0$ or to 
$2 \, \cosh^{-2} s$ on compact subsets. 
It is then easy to see that, some subsequence of $(\tilde{w}_n)_n$ 
converges to $w_\infty$, a nontrivial solution of one of the 
following equations~:
\[ 
\Delta_0 w_\infty  = 0, 
\] 
or   
\[ 
\Delta_0 w_\infty +\frac{2}{\cosh^2 (s+\bar{s})}w_\infty = 0, 
\qquad \mbox{for some} \quad\bar{s}\in {\mathbb R} , 
\] 
on $\R \times S^{1}$. In addition, $|w|\leq c\, 
e^{\tilde \mu s}$.  This is once again impossible. To see this, decompose
$w_\infty = \sum_j w_j (s) \, e^{i j \theta}$. By the choice of 
$\tilde \mu$ we get $w_0 = w_{\pm 1} =0$ since nontrivial solutions of 
the homogeneous problems on the eigenspaces $ j = 0, \pm 1$ decay at 
most like $\cosh^{-1}s$ at $\infty$, whereas $\tilde \mu  \in (-2,-1)$.
On the other hand, the restrictions of these two operators to the eigenspaces with 
$|j| \geq 2$ satisfy the maximum principle, and so all the remaining
components $w_j=0$. This is again a contradiction.

\medskip 

\noindent 
We have now ruled out all possibilities, and so the surfaces
${\mathfrak G}_\Sigma (p, \tau)$ are nondegenerate when $\tau$ is 
sufficiently small.

\section{Analysis of the forgetful map}
\label{anal-fm}
In the remainder of this paper we apply this gluing construction to 
study some aspects of the topology of $\calM_{g,k}$; this will be
accomplished using a natural mapping from this space into 
the more familiar and better understood Teichm\"uller space of closed 
Riemann surfaces of genus $g$ with an ordered $k$-tuple of points deleted.
We let $\calT_{g,k}^{{}}$ denote this latter space;  
any element has the form $\Sigb  -  \{p_1, \ldots, p_k\}$, where 
$\Sigb$ is a compact Riemann surface and the $p_j$ are distinct points on it. 

Any CMC surface $\Sig \in \calM_{g,k}$ is conformally equivalent to some 
element of this Teichm\"uller space. This allows us to define a 
forgetful map
\begin{equation}
\calF_{g,k}: \calM_{g,k} \longrightarrow \calT_{g,k}^{{}},
\label{eq:fm}
\end{equation}
which is given by forgetting the geometric (i.e. CMC) structure of
a surface and remembering only its conformal structure. For convenience,
we shall usually drop the subscripts $(g,k)$ from the notation
for this map. 

Our basic goal is to understand the image $\calI_{g,k} = \calI$
of this mapping. 
We first show that $\calF$ is real analytic. A recent 
result of Kusner \cite{Ku} states that (a slight modification of) 
$\calF$ is proper. Together, these results show that $\calI$ is a closed, 
subanalytic set. As such, it is stratified by real analytic submanifolds, 
and so we may 
define the codimension of $\calI$ as the codimension of its stratum of
maximal dimension. Using the preceding gluing construction,
we prove that $\calF$ is surjective when $g=0$ while for each
$g>0$ its image has codimension which is uniformly bounded 
(depending on $g$) as $k \to \infty$. Examined more carefully,
this argument also shows that $\calM_{g,k}$ detects much of the topology 
of this Teichm\"uller space. We conclude by summarizing the ramifications
of these results for the differential topological structure of the CMC 
moduli space. 

We stress that the basic ideas here very intuitive, granting the
main gluing theorem. Roughly speaking, the fact that $\calF_{0,k}$ is 
surjective when $g=0$ follows inductively from the fact that we can glue 
a half-Delaunay surface at any point $p$ of any fixed nondegenerate $\Sig 
\in \calM_{0,k'}$, $3 \leq k' \leq k$, without changing the conformal 
structure much away from $p$. The result is complicated when $g>0$ by 
the fact that there may be constraints on the image of $\calF$ which we do
not see directly; this is where the real analyticity enters, for it
implies that the image does lie in a well-behaved set. The uniform 
boundedness of the codimension of the image means essentially that
the only serious constraints on the image occur when $k$ is small,
and that when $k$ is large enough, the (conformal) location of the 
ends may be chosen freely. 

\subsection{Analyticity of $\calF$}

We show now that $\calF$ is analytic, not only when near smooth
points of $\calM_{g,k}$, but even as a map of (possibly singular) real 
analytic spaces. 

To state this result precisely, fix $\Sig \in \calM_{g,k}$ and let $g$ 
denote its induced metric. We parametrize a neighborhood of $\Sig$ in 
the space of all nearby surfaces in the usual way, using a finite 
dimensional family of deformations of $\Sig$ which preserve the CMC 
structure of the ends and then taking normal perturbations of these. 
Part of the main theorem in \cite{KMP} is that for any $\Sig \in
\calM_{g,k}$, there is a neighborhood in the CMC moduli space which lies 
in some open finite dimensional analytic submanifold $\calQ$ of this 
infinite dimensional space. When $\Sig$ is nondegenerate, then we may
assume that $\calQ$ is a real analytic coordinate chart in $\calM_{g,k}$, 
but in 
general, the CMC moduli space is the zero set of a real analytic function 
on $\calQ$. To be definite, we regard each element of 
$q$ as an embedding of a fixed surface $\Sig_0$ into $\R^3$; we let 
$q_0$ denote the base embedding, corresponding to the original CMC 
surface $\Sig$. From that construction, we may even assume that every 
$q \in \calQ$ is an analytic embedding, but this is not so important here. 

\begin{proposition} Suppose $\Sig \in \calM_{g,k}$, and let the
finite dimensional analytic manifold $\calQ$ be chosen as above.
Then the natural extension of the forgetful map $\calF$ which
assigns to any $q \in \calQ$ the element of $\calT_{g,k}^{}$
determined by $q(\Sig_0)$ is a real analytic mapping.
\label{pr:ra}
\end{proposition}

\noindent Since the CMC moduli space (locally) lies in $\calQ$, 
this result gives what is perhaps the most natural meaning to the
statement that $\calF$ is real analytic on $\calM_{g,k}$ 
near singular points of this moduli space.

\medskip

\noindent{\bf Proof:} 
Let $g_0$ denote the base metric on $\Sig_0$, i.e. $g_0 = q_0^*(\delta)$,
where $\delta$ is the Euclidean metric. Similarly, for $q \in \calQ$ we let 
$g_q = q^*(\delta)$. For each one of these metrics, there is a uniquely
determined conformal factor $e^{2\phi_q}$ such that $h_q = e^{2\phi_q}g_q$
is a hyperbolic metric of finite area. These hyperbolic metrics parametrize
the Teichm\"uller space $\calT_{g,k}^{}$, and so the theorem will
be proved if we show that the map $q \mapsto \phi_q$ is real analytic. 

For each end $E_j$ of $\Sig_0$, fix isothermal coordinates 
$(s_j,\theta_j)$. The model Delaunay surface $D_{\tau_j}$ for this
end has the metric
\[
g_{\tau_j} = (\tau_j^2 \, e^{2\sigma_{\tau_j}})(ds_j^2 + d\theta_j^2),
\]
and hence 
\[
g_q = (\tau_j^2 \, e^{2\sigma_{\tau_j}})\left(ds_j^2 + d\theta_j^2 + 
{\mathcal O}(e^{-\alpha s_j})\right)
\]
there, for some $\alpha > 0$. The Delaunay parameters $\tau_j$ and the
functions $\sigma_{\tau_j}$ depend real analytically on $q$. We can even 
modify this conformal factor further and thereby find a function $\mu_q$ 
which depends real analytically on $q$ such that
\[
g_{c,q} := e^{2\mu_q}g_q = ds_j^2 + d\theta_j^2 \qquad \mbox{on}\quad E_j.
\]
This is a metric with cylindrical ends. Accordingly, decompose 
the sought-after conformal factor $\phi_q$ as $\mu_q + \psi_q$.
We must prove that $\psi_q$ depends analytically on $q$. 

Letting $\Delta_q$ and $K_q$ denote the Laplace operator and Gauss
curvature function for $g_{c,q}$, then $\psi_q$ is the unique solution 
of the PDE 
\[ 
\Delta_{q}\psi_q - e^{2\psi_q} = K_{q}.
\]
Write $\psi_0$ for the unique solution when $q=q_0$. Note that 
$K_q = 0$ on all ends, and so $-\log s_j$ is a solution
along each $E_j$. 

We shall define natural function spaces $X$ and $Y$ below (these will 
be certain weighted Sobolev spaces), on which the mapping 
\[
\calN: \calQ \times X \longrightarrow Y,
\]
defined by 
\[
(q,\psi) \quad  \longmapsto \quad \Delta_{q}\psi - e^{2\psi} - K_{q}
\]
is locally surjective near $(q_0,\psi_0)$. To this end, observe that 
the differential in the second factor is
\[
L w := D_2\calN_{q_0,\psi_0}(w) = \Delta_{0}w - 2 e^{2\psi_0}w.
\]
In terms of the coordinates $(s_j,\theta_j)$ on $E_j$, 
\[
\left. L \right|_{E_j} = \del_{s_j}^2 + \del_{\theta_j}^2 - 2/{s_j}^2.
\]

It remains to choose the function spaces $X$ and $Y$ which have
the correct properties, in particular that $L:X \mapsto Y$ is surjective. 
While this is not too difficult, and equivalent theorems surely exist 
elsewhere in the literature, we sketch the proof for completeness.

As a first attempt, we can try to use a procedure similar to the
one outlined in \S 3.1; namely, if we define the weighted Sobolev
spaces $e^{ms} \, H^t(\Sig)$ (so $u$ is in this space if $u \in 
H^t_{\mathrm{loc}}$ and on each end $E_j$, $u = e^{ms_j}v$ where $v$ 
is in the ordinary Sobolev space $H^t(E_j)$), then it is straightforward 
to show that 
\begin{equation}
L: e^{ms} \, H^{t+2}(\Sig) \longrightarrow e^{ms} \, H^t(\Sig)
\label{eq:fl}
\end{equation}
is Fredholm whenever $m \notin {\mathbb Z}$. Furthermore, using
the maximum principle, which applies since the term of order
zero in $L$ is strictly negative, the mapping (\ref{eq:fl})
is injective when $m < 0$; by duality and elliptic
regularity, it is then surjective when $m > 0$, $m \neq 1,2,\ldots$. 
Again as in \S 3.1, there is a deficiency space consisting 
of the linear span of cutoffs of a special set of temperate solutions to 
$Lu = 0$ along each end. We can determine these temperate
solutions using separation of variables and see that they are linear 
combinations of the functions $s^2$ and $s^{-1}$. (For convenience, in 
many places in the remainder of this proof we drop references to the various
ends $E_j$ and also omit the subscript $j$.) Therefore, if $-1 < m < 0$,
then for any $f_1 \in e^{ms}H^t(\Sig)$ it is possible to find
a function 
\[
u_1 = a s^2 + b s^{-1} + \tilde u, \qquad \tilde u \in e^{ms}H^{t+2}(\Sig),
\]
where $a$ and $b$ are constants, and such that $Lu_1 = f_1$ on $\Sig$. 
The deficiency space $W$ consists of all linear combinations
of these solutions $s^2$ and $s^{-1}$ on all ends, and thus is 
$2k$-dimensional. We have shown that 
\begin{equation}
L: e^{ms}H^{t+2}(\Sig) \oplus W \longrightarrow e^{ms}H^t(\Sig)
\label{eq:F2}
\end{equation}
is surjective. A relative index calculation (cf.\ \cite{KMP}) shows
that the nullspace ${\mathcal B}$ of (\ref{eq:F2}) is $k$-dimensional. 

Unfortunately, this is not the end of the story because the nonlinear 
operator $\calN$ does not carry this domain space into $e^{ms}H^t(\Sig)$;
in fact, there is no evident way to use the geometric context to 
regularize this mapping. Therefore we proceed further.

\medskip

We first require Sobolev spaces with polynomial rather than exponential
weights. Thus for $t \in {\mathbb N}$ and $\nu \in {\R}$, let 
$H^t_\nu(\R^+ \times S^1)$ be the space of functions in $H^t_{loc}$ 
such that $s^{-1/2+k-\nu} \, \nabla^k u \in L^2 (\R^+ \times S^1)$, $k=0, 
\ldots, t$. We also define $H^t_\nu(\Sig)$ as the space of $H^t_{loc}(\Sig)$ 
functions which lie in $H^t_\nu(\R^+ \times S^1)$ on each end.

\medskip

By separation of variables, it is easy to produce a map
\[
G_0: H^t_{\nu-2}(\R^+ \times S^1) \longrightarrow H^{t+2}_{\nu}
(\R^+ \times S^1)
\]
for any $\nu \in \R$, $\nu \neq   -1, 2$, such that $(\del_s^2 + \del_\theta^2 - 2/s^2) \,  G_0 \,  
f_0 = f_0$. In other words, $G_0$ is a right inverse for $L$. Observe that we do not impose 
any boundary data, in particular $G_0$ is not unique. 
Using $G_0$ and cutoff functions $\chi_j$ on each ends, we can produce an operator 
\[
\tilde G_0: H^t_{\nu-2}(\Sig) \longrightarrow H^{t+2}_{\nu} (\Sig)
\]
 such that $u_0 : =  \tilde G_0 f_0$ has the property that $Lu_0$ 
vanishes outside a compact set of $\Sig$. 

Now, if $f \in H_{-2}^t(\Sig)$ and $t \geq 2$, $f$ is continuous.
Since the constant function $1$ is a subsolution for the 
operator $L$ -- in fact $L(1) = -2/s^2$ on each end -- we can
solve a sequence of equations $Lu_j = f$ in $\Sigma_j$  with $u_j = 0$ on the
boundaries $\del \Sigma_j$, where $\Sigma_j$ is a smooth exhaustion
of $\Sig$ by compact sets; the limit of this sequence is a {\it bounded}
solution $u$ of $Lu = f$. By earlier remarks, this is the only
bounded solution of this equation. If $f = 0$ on an end, then
we know that $u$ must be a linear combination of $s^2$ and 
$1/s$ and a term which exponentially decreases on that end, and
so by boundedness, $u = u_1 + v$ where $u_1 = a/s$ and $v \in
e^{-s} \, H^{t+2}$. 

Finally, let $\nu \in (-1, 0]$, $t \geq 2$ and define
\[
X = H_\nu^{t+2}(\Sig), \qquad \mbox{and}\qquad Y = H_{\nu-2}^t(\Sig).
\]
We claim that $L: X \rightarrow Y$ is an isomorphism and $\calN: 
\calQ \times X \rightarrow Y$ is smooth. To prove these, first
suppose that $f \in Y$. Let $u_0 : = \tilde G_0 f$, 
then $L(u-u_0) = \tilde{f}$ has compact support.  Next, there is also a unique 
bounded solution $\bar{u}$ of  $L\bar u = \tilde{f}$. As explained above,
 separating variables on the end, we see that 
$\tilde{v}: =  u_0 +\bar u$ is the sum of some multiple of $1/s$ and a
function in $H_\nu^{t+2}(\Sig) $.  Tracing through this procedure,
we have found a bounded map $G: Y \rightarrow X$ such that
$LG = I$. Since $L$ does not have any nullspace in $X$, this
map is an isomorphism. So far, we have only used the fact that 
$t \geq 2$ and $\nu \in (-1, 2)$.

We also note that $\calN$ is a real analytic mapping from
a neighbourhood of $0$ in $X$ to $Y$. This follows from the 
considerations above concerning the linear part $L$ of $\calN$, as 
well as the fact that
the nonlinear error term, which has the form 
$s^{-2} \, (e^{2w} - 1 - 2w)$ on each end.  This is where the 
restriction that $\nu \leq 0$ is required. 

In any case, we may now apply the real analytic version of the
implicit function theorem to get the existence of a real analytic 
mapping $\Psi: \calQ \to X$ such that $\calN(q,\Psi(q)) \equiv 0$. Since
$\psi_q = \Psi(q)$, we have proved the theorem.
\hfill $\Box$

\subsection{The image of $\calF$}

Recall that a connected component of $\calM_{g,k}$ is said to be 
nondegenerate if it contains an element $\Sig$ which is (analytically) 
nondegenerate in the sense of Definition~\ref{de:nond}. It is proved 
in \cite{MPPR} 
that $\calM_{g,k}$ contains a nondegenerate component for every $(g,k)$
with $g \geq 0$ and $k \geq 3$. The principal stratum in a nondegenerate 
component of $\calM_{g,k}$ has dimension $3k$. On the other hand, 
$\calT_{g,k}^{}$ is a real analytic manifold of dimension $6g-6 + 2k$. 
Therefore one might hope that $\calF$ is surjective, at least 
when $k$ is sufficiently large. 

We now give some results concerning the nature and size of the image 
$\calI$. These require some preliminary definitions.

There is a tautological bundle ${\mathbb V}_{g,k}$ over $\calT_{g,k}^{}$
defined by
\[
{\mathbb V}_{g,k} = \left\{([(\Sigb,p_1, \ldots,p_k)],p) \quad | \quad 
[(\Sigb,p_1, \ldots,p_k)] \in \calT_{g,k}^{},\ p \in
\Sigb  -  \{p_1, \ldots, p_k\} \right\}.
\]
This is the domain of a natural augmentation map
\[
\begin{aligned}
\mathcal A: {\mathbb V}_{g,k} \quad & \longrightarrow \quad 
\calT_{g,k+1}^{} \\
 ([(\Sigb,p_1, \ldots,p_k)],p) \quad & \longmapsto \quad
[(\Sigb,p_1, \ldots, p_k, p)].
\end{aligned}
\]
It is clear that $\mathcal A$ is an isomorphism.

Next, there is also a tautological bundle over $\calM_{g,k}$,
\[
{\mathbb U}_{g,k} = \{(\Sig,p) \quad | \quad  \Sig \in \calM_{g,k}, p \in \Sig \}.
\]
For any $(\Sig,p) \in {\mathbb U}_{g,k}$, let $(0,\tau^*(\Sig,p))$ 
denote the largest open interval such that if $0 < \tau < \tau^*(\Sig,p)$, 
then the gluing map which attaches a half-Delaunay end with Delaunay 
parameter $\tau$ to the point $p$ exists. It follows from the gluing
construction that $\tau^*(\Sig,p)$ is bounded away from zero on compact sets of $\Sigma$, provided
$\Sig$ lies in a compact set in a nondegenerate stratum of 
$\calM_{g,k}$. Let 
\[
\calW_{g,k} = \{(\Sig,p,\tau) \quad | \quad  (\Sig,p) \in {\mathbb U}_{g,k}, 0 < \tau < 
\tau^*(\Sig,p)\}.
\]
This is the natural domain of the gluing map
\[
{\mathfrak G} : \calW_{g,k} \longrightarrow \calM_{g,k+1}.
\]
This map is continuous, and in fact, real analytic. We shall often
omit the subscripts $(g,k)$ from these bundles when the meaning is clear.
Also, if $C$ is any subset either of $\calM_{g,k}$ or of 
${\mathbb U}_{g,k}$, then we let $\calW(C)$ denote the portion of 
$\calW$ lying over $C$; in particular $\calW(\Sig)$ denotes the 
natural domain of
the gluing map over a fixed surface $\Sig$. We write the 
CMC surface ${\mathfrak G}_\Sigma (p,\tau)$ as $\Sig_{p,\tau}$. Finally, note that
\begin{equation}
\lim_{\tau \to 0} \calF(\Sig_{p,\tau}) = {\mathcal A}(\calF(\Sig),p).
\label{eq:aug}
\end{equation}

\begin{theorem} Suppose that $g=0$ and $k \geq 3$. Then there is
a nondegenerate component of $\calM_{0,k}$ on which $\calF$ is 
surjective.
\end{theorem}

\noindent{\bf Proof:} We prove this by induction on $k$. When
$k=3$, then according to \cite{GKS}, $\calM_{0,3}$ is homeomorphic 
to a $3$-ball, hence in particular is connected; by \cite{MP} it contains 
a nondegenerate element. In other words, $\calM_{0,3}$ contains a 
single component, and this component is nondegenerate. On the other hand, 
$\calT_{0,3}^{}$ consists of a single point. Hence $\calF_{0,3}$ 
is obviously surjective.

Now suppose that $k \geq 3$, and $C_k^0$ is a nondegenerate stratum 
in some component $C_k \subset \calM_{0,k}$ such that 
$\calF: C_k^0 \to \calT_{0,k}^{}$ is surjective. Choose any
point $[(S^2,p_1,\ldots,p_k,p_{k+1})] \in \calT_{0,k+1}^{}$
(note that $\Sigb$ must be $S^2$ when $g=0$), and let
$([(S^2,p_1,\ldots,p_k)],p_{k+1})$ be the lift of this point to
${\mathbb V}_{0,k}$. Write $p$ for $p_{k+1}$ for simplicity.  

By assumption, there is an element $\Sig \in C_k^0$ 
such that $\calF(\Sig) = [(S^2,p_1,\ldots,p_k)]$. 
Let $\mathcal B$ be some neighbourhood of $(\Sig,p)$ in
${\mathbb U}_{0,k}$. Then for any $(\Sig',q,\tau) \in \mathcal 
B \times (0,\eta)$, we obtain elements $\Sig'_{q,\tau} \in 
\calM_{0,k+1}$ and $(\calF(\Sig'),q)\in {\mathbb V}_{0,k}$. 
The theorem will be proved if we show that 
\[
{\mathcal A}(\calF(\Sig),p) \in (\calF \circ {\mathfrak G} )
({\mathcal B},\eta),
\]
when $\eta$ is small enough. But this is clear from 
(\ref{eq:aug}) using a straightforward degree theory argument.  

This proof also shows that the nondegenerate stratum $C_{k+1}^0$
is obtained inductively by gluing half-Delaunay surfaces with
very small necks located at arbitrary points $p \in \Sig$ for all 
surfaces $\Sig \in C_k^0$. \hfill $\Box$

\medskip

As remarked above, when $g > 0$ we no longer expect $\calF_{g,k}$
to be surjective. We shall show instead that its codimension
does not become unbounded as $k \to \infty$. 
\begin{theorem} Suppose that $g \geq 1$ and that $C_{k_0}$ is a nondegenerate
component of $\calM_{g,k_0}$, with nondegenerate stratum $C_{k_0}^0$.
Then for each $k = k_0, k_0 + 1, \ldots$ there is a nondegenerate component 
$C_k \subset \calM_{g,k}$ such that codimension of the image 
$\calI_k = \calF(C_k)$ is bounded as $k\to \infty$.
\end{theorem}

\noindent{\bf Proof:} This is also proved by an inductive procedure.
Let us suppose that for some $k \geq k_0$, there is a nondegenerate
element $\Sig \in \calM_{g,k}$; suppose also that $\mbox{rank}\,
\left.D\calF\right|_\Sig \equiv r$ is maximal amongst all such elements. 
Let $C_k$ be the component of $\calM_{g,k}$ containing $\Sig$
and $\calI_k$ its image by $\calF$ in $\calT_{g,k}^{}$;
The dimension of the stratum of $\calI_k$ through $\calF(\Sig)$ is $r$,
and so $d_k = \mbox{codim}\,(\calI_k) \leq 6g-6 + 2k - r$.
Let us furthermore choose an $r$-dimensional analytic submanifold 
$\calS_k$ through $\Sig$ such that the restriction of $\calF$ to it 
is an analytic diffeomorphism onto its image. Let ${\mathbb U}(\calS_k)$
be the portion of the bundle ${\mathcal U}$ lying over $\calS_k$.

We now consider, for some small $\eta>0$, the restriction of the gluing map
\[
{\mathfrak G}_\eta: {\mathbb U}(\calS_k) \longrightarrow \calM_{g,k+1}.
\]
We first claim that the dimension of the image of ${\mathfrak G}_\eta$ is $r+2$.
This is straightforward, since the differential of ${\mathfrak G}_\eta$ in
the directions of the fibres of ${\mathbb U}(\calS_k)$, i.e. letting
$p$ vary and $\Sig' \in \calS$ remain fixed, is injective. 
This produces an $(r+2)$-dimensional submanifold $\calS_{k+1}$
in $\calM_{g,k+1}$, consisting entirely of nondegenerate points.
Using (\ref{eq:aug}) again, when $\eta$ is small enough,
$\calF(\calS_{k+1})$ is $(r+2)$-dimensional. However, since
$\dim(\calT_{g,k+1}^{}) - \dim(\calT_{g,k}^{}) = 2$,
we see that the codimension of $\calF(\calS_{k+1})$
is again $6g-6+2k-r$, and so the codimension of the image
of the component of $\calM_{g,k+1}$ containing $\calS_{k+1}$
is bounded by this same number. This proves the theorem.
\hfill $\Box$

\medskip

Regarding the global structure of the image $\calI_{g,k}$ of 
$\calM_{g,k}$ by $\calF$, we quote a recent nice result of
Kusner. To state it, recall that the necksize parameter $\tau$ of
an end $E$ of $\Sig \in \calM_{g,k}$ is the Delaunay parameter
of the Delaunay surface to which this end is asymptotic. This
value may be determined using the force integral. 

\begin{proposition} (Kusner \cite{Ku}) Let $\Sig_\ell$ be a sequence
of elements in $\calM_{g,k}$ and suppose that the necksize 
parameters $\tau_j^{(\ell)}$ of the ends $E_j$ of $\Sig_\ell$
are all bounded below by some $\eta > 0$. Then either the
conformal structures $\calF(\Sig_\ell)$ diverge in $\calT_{g,k}^{}$
or else the surfaces $\Sig_\ell$ converge, up to rigid motion,
to some limiting surface $\Sig_\infty \in \calM_{g,k}$. 
\label{pr:Kusner}
\end{proposition}

\noindent Kusner's result is somewhat more general, and he phrases
it in terms of properness of $\calF$. 

\medskip

Combining Propositions~\ref{pr:ra} and \ref{pr:Kusner}, we obtain
\begin{proposition}
The image $\calI_{g,k}$ of the forgetful map $\calF$, restricted to
any (not necessarily nondegenerate) component in $\calT_{g,k}^{}$
is a closed subanalytic set.
\end{proposition}

We conclude this subsection with a discussion of the fundamental group 
of $\calM_{g,k}$. While we do not determine this group precisely,
we examine the homomorphism
\[
\calF_*: \pi_1(\calM_{g,k}) \longrightarrow
\pi_1(\calT_{g,k}^{}).
\]
The group on the right here is well understood and rather complicated.
We show that the image of $\calF_*$ is a fairly large subgroup. 

\medskip

We first review some facts about the group $\pi_1(\calT_{g,k}^{})$,
refering to \cite{Bi} for more details. There is a subsidiary
forgetful map
\[
F': \calT_{g,k}^{} \longrightarrow \calC(g,k).
\]
The space on the right here is the Teichm\"uller space
of conformal structures on a compact surface of genus $g$ (i.e.
the classical Teichm\"uller space) and the configuration space
of $k$ points on a compact surface of genus $g$. To define it,
recall that we may identify an element of $\calT_{g,k}^{}$
with a hyperbolic metric on the compact surface $\Sigb$ along
with an ordered $k$-tuple of distinct points $(p_1, \ldots, p_k)$ 
on $\Sigb$ (rather than finite area complete hyperbolic metrics
on $\Sigb - \{p_1, \ldots, p_k\}$). Then $F'$ 
 is defined by forgetting the conformal structure.
It is well-known that 
\[
\calF'_*: \pi_1(\calT_{g,k}^{}) \longrightarrow \pi_1(\calC(g,k))
\]
is an isomorphism. We write $\calF' = F' \circ \calF$. 

\medskip

The space $\calC(g,k)$ has a rather interesting topology. Its 
fundamental group is known as the pure braid group of the surface 
of genus $g$ on $k$ braids, and we denote it by $B(g,k)$. It is a 
finitely generated group; the loops $\Gamma_{ij}$, $1 \leq i,j 
\leq k$, $i \neq j$, corresponding to the point $p_j$ traversing a 
small loop winding around the point $p_i$ once, with all other $p_\ell$ 
fixed, comprise a generating set. 

\begin{theorem} When $k \geq 3$, the map 
\[
\calF'_*: \pi_1(\calM_{0,k}) \longrightarrow \pi_1(\calC(0,k))
\]
is an epimorphism.
\label{pi1a}
\end{theorem}

\begin{theorem} Suppose $g \geq 1$ and $\calM_{g,k_0}$ contains
a nondegenerate component. Then for any $k \geq k_0$,
$\calM_{g,k}$ contains a nondegenerate component $C_k$ such
that the image of the homomorphism
\[
\calF'_*: \pi_1(C_k) \longrightarrow \pi_1(\calC(g,k))
\]
contains the subgroup of $B(g,k)$ generated by the collection of
loops $\Gamma_{ij}$, $j > k_0$.
\label{pi1b}
\end{theorem}

The proofs of these two theorems are nearly identical. They
rely only on the simple observations that any of the 
loops $\Gamma_{ij}$ are in the image of $\calF'_*$ when $g=0$, while
when $g>0$, at least those loops with $j>k_0$ are in the image. 

\subsection{The structure of the CMC moduli space}

We conclude this paper by describing informally what we know 
at this point about the CMC moduli spaces $\calM_{g,k}$. 

As before, we let $\calI_{g,k}$ denote the image of $\calM_{g,k}$
under the forgetful map $\calF$. Then it is tempting to think of
\[
\calF: \calM_{g,k} \longrightarrow \calI_{g,k} \subset \calT_{g,k}^{}
\]
as a sort of singular fibration. We have shown that all spaces here
are real analytic or subanalytic, hence stratified, and $\calF$ is 
a real analytic mapping. The image $\calI_{g,k}$ detects at least
some fairly large portion of the fundamental group of $\calT_{g,k}^{}$
when $g \geq 1$ and $k$ is large; it detects all of it when $g=0$
and $k \geq 3$. If $\calM_{g,k}(\eta)$ denotes the subset of surfaces
$\Sig \in \calM_{g,k}$ with necksizes of all ends of $\Sig$ no smaller
than $\eta$, then by Kusner's theorem, the restriction of $\calF$
to this subset is proper. As already noted, one way to interpret
this is that if $\Sig_j$ is a divergent sequence of surfaces in 
$\calM_{g,k}$ with no end necksizes tending to zero, then
necessarily the conformal structures $\calF(\Sig_j)$ must
be degenerating. This behaviour does indeed occur; for example,
the connected sum construction of \cite{MPP}, cf.\ also \cite{MPPR}, 
shows that it 
is possible to construct sequences of surfaces with no end
necksizes tending to zero, but with some interior necks pinching
off. This corresponds to degeneration in $\calT_{g,k}^{}$,
and these examples exist even when $g=0$. 

\medskip

We conclude with a number of open questions~:

\begin{itemize}

\item Is $\calM_{g,k}$ connected? The only case where this is 
understood (and known to be true) is when $(g,k) = (0,3)$ by \cite{GKS}.

\item Do the fibres of $\calF$ ever have nontrivial topology; for example,
do they ever contain homotopically nontrivial 
loops? This does not occur in $\calM_{0,3}$, and if the fibres
are always contractible, then the image $\calI_{g,k}$ of
any component of $\calM_{g,k}$ would be a retract of that
component. 

\item Is $\calF_{g,k}$ ever surjective when $g>0$? A heuristic argument
against this might be made by considering configurations
$(\Sigb,p_1, \ldots, p_k)$ where all of the points $p_j$ are
contained in a small neighbourhood (for example, a small ball
relative to the conformally equivalent flat or hyperbolic
metric on $\Sigb$); it seems likely that the CMC balancing
formul\ae\ would rule out CMC realizations of such configurations.

\item Construct, or prove the existence of, a degenerate 
CMC surface $\Sig$ in some $\calM_{g,k}$. Although the possibility 
of their existence adds significant complications throughout the theory, 
to date none are known to exist.

\item In a related direction, it seems quite likely that 
every element in $\calM_{0,3}$ is nondegenerate, and it would
be very useful to know whether this is true. One motivation is
that the (otherwise) very explicit geometric knowledge about these 
surfaces makes them ideally suited as building blocks in gluing
constructions, but to use them in this way requires knowing
that they are nondegenerate.

\item We have not discussed the geometric structure on these
CMC moduli spaces. Along these lines, it is proved in \cite{KMP}
that on the infinitesimal level, each $\calM_{g,k}$ has the
structure of a Lagrangian submanifold of a larger symplectic
submanifold. More precisely, the tangent space of $\calM_{g,k}$
at any nondegenerate point is a Lagrangian subspace of a natural
symplectic vector space. It is not too difficult to make this
global picture more precise, but the more compelling question
is~: what can be done with it? Some simple examples and other
evidence point to the possibility that there is a tautological 
one-form on some large subset of this ambient symplectic
manifold which becomes singular on the subvariety of surfaces
with at least one (asymptotically) cylindrical end. The periods
of this one-form appear to have direct geometric meaning.
Is there anything else which can be done with this Lagrangian
structure?
\end{itemize}


\begin{thebibliography}{99}

\bibitem{Bi} J. Birman,
{\em Braids, Links and Mapping Class Groups}, 
Annals of Mathematical Studies {\bf 82}
Princeton University Press, Princeton (1975).

\bibitem{FP} S. Fakhi and F. Pacard,
{\em Existence result for minimal hypersurfaces with a prescribed finite 
number of planar ends,}
Manuscripta Math. {\bf 103} (2000), no. 4, 465--512.

\bibitem{GKS} K. Grosse-Brauckmann, R. Kusner and J. Sullivan,
{\em Triunduloids: embedded constant mean curvatures surfaces with three ends 
and genus zero,} Preprint. Math.DG/0102183.

\bibitem{Ka} N. Kapouleas,
{\em Complete constant mean curvature surfaces in Euclidean three space,}
Ann. of Math. (2) {\bf 131} (1990), 239--330.

\bibitem{Ku} R. Kusner,
{\em Conformal structure of embedded CMC surfaces}, This volume.

\bibitem{KKS} N. Korevaar, R. Kusner and B. Solomon,
{\em The structure of complete embedded surfaces with constant 
mean curvature},
J. Differential Geometry 30 (1989) 465--503

\bibitem{KMP} R. Kusner, R. Mazzeo and D. Pollack,
{\em The moduli space of complete embedded constant mean curvature
surfaces,}
Geom. Funct. Anal. {\bf 6} (1996) 120--137.

\bibitem{MP} R. Mazzeo and F. Pacard,
{\em Constant mean curvature surfaces with Delaunay ends,}
Comm. Anal. Geom. {\bf 9} No. 1  (2001) 169--237.

\bibitem{MP2} R. Mazzeo and F. Pacard, {\em Bifurcating nodoids},
To appear. 

\bibitem{MPP} R. Mazzeo, F. Pacard and D. Pollack,
{\em Connected Sums of constant mean curvature surfaces in Euclidean 3
space,} J. Reine Angew. Math. {\bf 536} (2001), 115--165.
 
\bibitem{MPPR} R. Mazzeo, F. Pacard, D. Pollack and J. Ratzkin, 
In Preparation. 

\bibitem{MPU} R. Mazzeo, D. Pollack and K. Uhlenbeck, 
{\em Moduli spaces of singular Yamabe metrics,}
J. Amer. Math. Soc. {\bf 9} (1996), no. 2, 303--344. 

\bibitem {Rat} J. Ratzkin, {\em An end-to-end gluing construction for surfaces
of constant mean curvature}, PhD Thesis, University of Washington (2001).

\bibitem{Mee} W. Meeks III,
{\em The topology and geometry of embedded surfaces of constant mean
curvature,} 
J. Differential Geom. {\bf 27} (1988), no. 3, 539--552.

\end{thebibliography}
\end{document}